\documentclass[reqno,11pt]{article}
\usepackage{bbm}
\usepackage{mathrsfs}
\usepackage{amssymb}
\usepackage[usenames,dvipsnames]{xcolor}
\usepackage{amsthm}
\usepackage{amsmath}
\usepackage{amsfonts}
\usepackage{enumerate}
\usepackage{txfonts}
\usepackage{bm}

\usepackage{graphicx}

\numberwithin{equation}{section}

\newtheorem{theorem}{Theorem}[section]
\newtheorem{lemma}[theorem]{Lemma}
\newtheorem{corollary}[theorem]{Corollary}
\newtheorem{proposition}[theorem]{Proposition}
\newtheorem{problem}[theorem]{Problem}
\newtheorem{remark}[theorem]{Remark}

\theoremstyle{definition}
\newtheorem{definition}[theorem]{Definition}
\newtheorem{assumption}[theorem]{Assumption}
\newtheorem{example}[theorem]{Example}

\begin{document}

\def\be{\begin{eqnarray}}
\def\ee{\end{eqnarray}}
\def\p{\partial}
\def\no{\nonumber}
\def\e{\epsilon}
\def\de{\delta}
\def\De{\Delta}
\def\om{\omega}
\def\Om{\Omega}
\def\f{\frac}
\def\th{\theta}
\def\la{\lambda}
\def\lab{\label}
\def\b{\bigg}
\def\var{\varphi}
\def\na{\nabla}
\def\ka{\kappa}
\def\al{\alpha}
\def\La{\Lambda}
\def\ga{\gamma}
\def\Ga{\Gamma}
\def\ti{\tilde}
\def\wti{\widetilde}
\def\wh{\widehat}
\def\ol{\overline}
\def\ul{\underline}
\def\Th{\Theta}
\def\si{\sigma}
\def\Si{\Sigma}
\def\oo{\infty}
\def\q{\quad}
\def\z{\zeta}
\def\co{\coloneqq}
\def\eqq{\eqqcolon}
\def\bt{\begin{theorem}}
\def\et{\end{theorem}}
\def\bc{\begin{corollary}}
\def\ec{\end{corollary}}
\def\bl{\begin{lemma}}
\def\el{\end{lemma}}
\def\bp{\begin{proposition}}
\def\ep{\end{proposition}}
\def\br{\begin{remark}}
\def\er{\end{remark}}
\def\bd{\begin{definition}}
\def\ed{\end{definition}}
\def\bpf{\begin{proof}}
\def\epf{\end{proof}}
\def\bex{\begin{example}}
\def\eex{\end{example}}
\def\bq{\begin{question}}
\def\eq{\end{question}}
\def\bas{\begin{assumption}}
\def\eas{\end{assumption}}
\def\ber{\begin{exercise}}
\def\eer{\end{exercise}}
\def\mb{\mathbb}
\def\mbR{\mb{R}}
\def\mbZ{\mb{Z}}
\def\mc{\mathcal}
\def\mcS{\mc{S}}
\def\ms{\mathscr}
\def\lan{\langle}
\def\ran{\rangle}
\def\lb{\llbracket}
\def\rb{\rrbracket}

\def\ol{\overline}
\def\der{\partial}
\def\tPsi{\tilde{\Psi}}
\def\tphi{\tilde{\phi}}
\def\tpsi{\tilde{\psi}}
\def\alp{\alpha}
\def\mcl{\mathcal}
\def\tx{\text}
\def\R{\mathbb{R}}
\def \vphi{\varphi}

\title{3-D axisymmetric subsonic flows with nonzero swirl for the compressible Euler-Poisson system}

\author{Myoungjean Bae\thanks{Department of Mathematics, POSTECH. Pohang, Gyungbuk, 37673, Republic of Korea ;
         Korea Institute for Advanced Study,
85 Hoegiro, Dongdaemun-gu,
Seoul 02455,
Republic of Korea
. Email: mjbae@postech.ac.kr.}\and {Shangkun Weng\thanks{Pohang Mathematics Institute, POSTECH. Pohang, Gyungbuk, 37673, Republic of Korea. Email: skwengmath@gmail.com.}}}

\pagestyle{myheadings} \markboth{3-D axisymmetric Euler-Poisson system with nonzero swirl}{3-D axisymmetric Euler-Poisson system with nonzero swirl}\maketitle

\begin{abstract}
  We address the structural stability of 3-D axisymmetric subsonic flows with nonzero swirl for the steady compressible Euler-Poisson system in a cylinder supplemented with non small boundary data.
A special Helmholtz decomposition of the velocity field is introduced for 3-D axisymmetric flow with a nonzero swirl(=angular momentum density) component.
 With the newly introduced decomposition, a quasilinear elliptic system of second order is derived from the elliptic modes in Euler-Poisson system for subsonic flows.
Due to the nonzero swirl, the main difficulties lie in the solvability of a singular elliptic equation which concerns the angular component of the vorticity in its cylindrical representation, and in analysis of streamlines near the axis $r=0$.
\end{abstract}

\begin{center}
\begin{minipage}{5.5in}
\end{minipage}
\end{center}

\section{Introduction and main results}

The steady Euler-Poisson system
\be\label{b11}
\begin{cases}
\text{div }(\rho {\bf u})=0,\\
\text{div }(\rho {\bf u}\otimes {\bf u}+ p I_n) =\rho \nabla \Phi,\\
\text{div }(\rho \mc{E} {\bf u}+ p{\bf u}) = \rho \bm{u}\cdot \nabla \Phi,\\
{\Delta_{{\rm x}}} \Phi= \rho -b({\rm x}).
\end{cases}
\ee
is a hydrodynamical model of semiconductor devices or plasmas, describing local behaviors of the electron density $\rho$, the macroscopic particle velocity ${\bf u}$, and the total energy $\mc{E}=|{\bf u}|^2/2+ e$, where $e$ is the internal energy. The first equation, which is also called as the continuity equation, expresses the conservation of electrons, the second equations express the conservation of momentum, where $\rho \na \Phi$ is the Coulomb force of electron particles. The third equation expresses the conservation of energy, and the last Poisson equation expresses the local change of the electric potential $\Phi$ due to the the volumetric charge density. The function $b({\rm x})>0$ is the prescribed density of fixed, positively charged background ions. Physically, by solving the Euler-Poisson equations in predetermined macroscopic device region with the relevant boundary conditions, we get the electric distribution or electric current in any proper cross sections.

To close the system \eqref{b11}, we introduce the equation of state
\be\lab{b12}
p=p(\rho,e)=(\ga-1) \rho e ,
\ee
where $\ga>1$ is called the {\emph{adiabatic constant}}. In terms of the entropy $S$, one also has
\be\lab{b14}
p(\rho, S)= A \text{exp }\b(\f{S}{c_v}\b) \,\rho^{\ga},
\ee
where $A$ and $c_v$ are positive constants. For more details about the physical background of the semiconductor device or models, one may refer to \cite{markowich1,markowich2,mrs}.

Define Bernoulli's function $\mc{B}$ by
\be\lab{b16}
\mc{B}= \f{|{\bf u}|^2}{2} + e + \f{p}{\rho} =\f{|{\bf u}|^2}{2} + \f{A \ga}{\ga-1} \textit{exp }\b(\f{S}{c_v}\b) \,\rho^{\ga-1}.
\ee
Then, the system (\ref{b11}) can be rewritten as
\be\lab{b15}
\begin{cases}
\text{div }(\rho {\bf u})=0,\\
\text{div }(\rho {\bf u}\otimes {\bf u}+ p I_n) =\rho \nabla \Phi,\\
\text{div }(\rho {\bf u}\mc{B})=\rho {\bf u} \cdot \nabla\Phi,\\
{\Delta_{{\rm x}}}  \Phi= \rho -b({\bf x}).
\end{cases}
\ee
The system (\ref{b15}) is a hyperbolic-elliptic coupled system, and behaves quite differently in subsonic states($|{\bf u}|<\sqrt{\p_{\rho}p(\rho, S)}$) and supersonic states($|{\bf u}|>\sqrt{\p_{\rho}p(\rho, S)}$), respectively. The goal of this work is to prove the structural stability of three dimensional axially symmetric subsonic flows with nonzero  swirl(=nonzero angular momentum) to the system \eqref{b15} for non small boundary data.

The existence and the uniqueness of subsonic flows to Euler-Poisson system were proved in \cite{abjm,amps, bdx,bdx14,bdx15,dm1, dm2,  markowich2, weng14, yeh}. In \cite{dm1, dm2}, the unique existence of subsonic flows for Euler-Poisson system is proved for small data. Subsonic flows with small current flux were studied in \cite{abjm,amps,  markowich2, yeh}. The structural stability of subsonic flows for multidimensional potential flow and two dimensional flow with nonzero vorticity was proved in \cite{bdx,bdx14,bdx15}, where no smallness of data was assumed. In \cite{weng14}, the unique existence of three dimensional subsonic flows with nonzero vorticity was proved. It used the Bernoulli's law to provide a new formulation of Euler-Poisson equations by reducing the dimension of the velocity, this idea is originally from \cite{weng15}. Although the method in \cite{weng14} works for the 3-D non-isentropic Euler-Poisson system, there are some smallness requirements on the background solutions.

The new feature of this work is that we construct three dimensional subsonic flows with nonzero vorticity, and that no smallness of data is required. In \cite{bdx}, it is found that a special structure of potential flow model of Euler-Poisson system yields the structural stability of multidimensional subsonic solutions without assumption of smallness of data. This result is extended to the case of two dimensional flow with nonzero vorticity  through a two dimensional Helmholtz decomposition ${\bf u}= \na \var+ \na^{\perp} \psi$ in \cite{bdx14}. In this paper, we introduce a Helmholtz decomposition for three dimensional subsonic flows in the form of
\be\no
{\bf u}=\nabla \varphi+{\rm curl} {\bf V}\quad\text{with}\,\,{\bf V}= h{\bf e}_r+\psi{\bf e}_{\theta},
\ee
where $\varphi, h, \psi$ are functions of $(x,r)$ for {$r=\sqrt{x_2^2+x_3^2}$}. With using this decomposition, we investigate axisymmetric subsonic flows with nonzero vorticity.  In particular, the function $\psi$ concerns the swirl (=angular momentum density). There are many other studies of axially symmetric smooth subsonic solutions to the steady compressible Euler equations \cite{acf,dd,lxy11,xx}. To our best knowledge, all the previous studies on the steady axially symmetric flows  mostly concern the case of zero swirl component.
We expect not only the result of this work contributes to understand a stabilizing or in-stabilizing effect of  vorticity to three dimensional subsonic flows of Euler-Poisson system, but also the new Helmholtz decomposition introduced in this work may open a new approach to investigate multidimensional transonic shock solutions to Euler-Poisson system or even transonic shock solutions to Euler system, which were previously studied in \cite{bf,chenf,cy, dwx, gamba, gm, lxy13, lxy11, lx, lrxx} and in the references therein.

\medskip

In the cylindrical coordinates $(x, r, \theta)$ satisfying
\be\no
(x_1, x_2, x_3)=(x, r\cos\theta, r\sin\theta),
\ee
for ${\rm x}=(x_1, x_2, x_3)\in \mathbb{R}^3$,
any function $f({\bf x})$ can be represented as $f({\bf x})=f(x,r,\theta)$, and a vector-valued function ${\bf u}({\rm x})$ can be represented as ${\bf u}({\rm x})=u_x(x,r,\th)\,{\bf e}_x+ u_r(x,r,\th)\,{\bf e}_r+ u_{\th}(x,r,\th)\,{\bf e}_{\th}$,
where
\be\no
{\bf e}_x=(1,0,0),\quad {\bf e}_r=(0,\cos\theta, \sin\theta),\quad {\bf e}_{\theta}=(0,-\sin\theta,\cos\theta).
\ee
We say that a function $f({\rm x})$ is {\it axially symmetric} if its value is independent of $\th$ and that a vector-valued function ${\bf h}= (h_x, h_r, h_{\th})$ is {\it axially symmetric} if each of functions $h_x({\bf x}), h_{r}({\bf x})$ and $h_{\th}({\bf x})$ is axially symmetric.

Assume that the smooth solution $(\rho, {\bf u}, S, \Phi)$ is axially symmetric, i.e.
\be\no
\rho({\rm x})&=&\rho(x,r),\q S({\rm x})=S(x,r),\q \Phi({\rm x})=\Phi(x,r),\\\no
{\bf u}({\rm x})&=& u_x(x,r) {\bf e}_{x}+ u_{r}(x,r) {\bf e}_{r}+ u_{\theta}(x,r) {\bf e}_{\theta},
\ee
then (\ref{b15}) can be simplified as
\be\lab{b17}
\begin{cases}
\p_x (\rho u_x) + \p_r (\rho u_r) + \f{\rho u_r}{r}=0,\\
\rho (u_x\p_x + u_r \p_r ) u_x + \p_x p = \rho \p_x \Phi,\\
\rho (u_x\p_x + u_r \p_r ) u_r - \f{\rho u_{\th}^2}{r} + \p_r p= \rho \p_r \Phi,\\
\rho (u_x\p_x + u_r \p_r ) u_{\th}+ \f{\rho u_r u_{\th}}{r}= 0,\\
\rho (u_x\p_x + u_r \p_r ) S=0,\\
\Delta \Phi = \rho -b,
\end{cases}
\ee
Hereafter we assume that the function $b$ is axially symmetric. Define
\be\no
\La(x,r)\co r u_{\th}(x,r).
\ee
$\La(x,r)$ represents the angular momentum density, and it is derived from \eqref{b17} that
\be\lab{b18}
\rho (u_x\p_x + u_r \p_r) \La =0.
\ee

Given a constant $L>0$, we fix a three dimensional axially symmetric nozzle of the length L by
\be\no
\mc{N}\co \{{\bf x}=(x_1,x_2, x_3)\in \mbR^3: 0< x_1< L, \,\,{x_2^2+x_3^2}< 1\}.
\ee
The entrance $\Ga_0$, exit $\Ga_L$ and the wall $\Ga_w$ of the nozzle $\mc{N}$ are defined as
\be\no
&&\Ga_0= \{0\}\times \{(x_2,x_3): x_2^2+x_3^2\leq 1\},\q \Ga_L= \{L\}\times \{(x_2,x_3): x_2^2+x_3^2\leq 1\},\\\no
&&\Ga_w= (0,L)\times \{(x_2,x_3): x_2^2+x_3^2= 1\}.
\ee
Fix a constant $b_0>0$. We first  compute one dimensional solutions $(\rho, {\bf u}, p, \Phi)$ of (\ref{b15}) in $\mc{N}$ with $b=b_0, u_2=u_3=\Phi_{x_2}=\Phi_{x_3}=0$, $\rho>0$, and $u_1>0$.
Set $E\co \Phi_{x_1}$. Then \eqref{b15} is reduced to the following ODE system for $(\rho, u_1, S, E)$
\be\lab{background1}\begin{array}{ll}
(\rho u_1)'=0,\\
S'=0,\\
\mc{B}'= E,\\
E'=\rho-b_0,
\end{array}\ee
where $'$ denotes the derivative with respect to $x_1$. Then $\rho u_1= J_0$ and $S=S_0$ for constants $J_0>0$ and $S_0>0$ to be determined by the entrance data. Therefore we can further reduce (\ref{background1}) to the following ODE system for $(\rho, E)$:
\be\lab{background2} \begin{array}{ll}
\rho'=\f{\rho E}{\ga A \textit{exp}\b(\f{S_0}{c_v}\b)\rho^{\ga-1}-\f{J_0^2}{\rho^2}},\\
E'=\rho-b_0.
\end{array}
\ee
For a detailed information on various types of solutions to \eqref{background2}, one can refer to \cite{lx}.

\bp\lab{background}
{\it Fix $b_0>0$. Given constants $J_0>0, S_0>0, \rho_0>\rho_c\co\b(\f{J_0^2}{\ga A\textit{exp}(\f{S_0}{c_v})}\b)^{\f{1}{\ga+1}}$ and $E_0$, there exist positive constants $L, \rho_{\sharp}, \rho^{\sharp}$, and $\nu_0$ such that the initial value problem (\ref{background2}) with
\be\lab{background3}
(\rho, E)(0)=(\rho_0, E_0)
\ee
has a unique smooth solution $(\rho(x), E(x))$ on the interval $[0, L]$ satisfying
\be\lab{background4}
\rho_c<\rho_{\sharp}\leq \rho \leq \rho^{\sharp}\q \textit{and}\q
\ga A\textit{exp}\b(\f{S_0}{c_v}\b)\, \rho^{\ga-1}-\f{J_0^2}{\rho^2}\geq \nu_0\q \textit{on }\q [0, L].
\ee
}
\ep
A proof of this proposition can be easily given by adjusting the argument in \cite[Section 1.3]{bdx}. In \cite{bdx}. it is shown that for any given length $L$ of  the nozzle, there exists a nonempty set $\mathfrak{P}$ of the entrance data $(\rho_0, E_0)$ so that a subsonic solution to \eqref{background2} with \eqref{background3} for each $(\rho_0, E_0)\in \mathfrak{P}$ uniquely exists.

\medskip

For fixed positive constants  $(b_0, J_0, S_0, \rho_0)$ with $\rho_0>\rho_c$, and a fixed constant $E_0$(which is not necessarily positive), let $(\bar{\rho}, \bar{E})$ be the unique smooth solution to the initial value problem of \eqref{background2} and \eqref{background3} on the interval $[0, L]$. Then we set $\bar{u}_1:=\frac{J_0}{\bar{\rho}}$ and $\bar{\bf u}=(\bar u_1, 0,0)$.

We define $\Phi_0({\rm x})$ and $\var_0({\rm x})$ by
\be\lab{background5}
\Phi_0({\rm x})= \int_0^{x_1} \bar E(t) \,dt\q \textit{and}\q \var_0({\rm x})= \int_0^{x_1} \bar u_1(t)\,dt\q \textit{for }{\rm x}=(x_1,x_2,x_3)\in\mc{N}.
\ee

We call $\bar{U}:=(\bar{\rho}, \bar{\bf{u}}, S_0, \Phi_0)$ the \emph{background solution} to \eqref{b15} in $\mc{N}$ \emph{associated with the entrance data $(b_0, J_0, S_0, \rho_0, E_0)$}. For the background solution, set
\be\lab{def-B0}
\mc{B}_0:=\frac{J_0^2}{2\rho_0^2}+\frac{\ga A}{\ga-1}{\exp}\left(\frac{S_0}{C_{\nu}}\right)\rho_0^{\ga-1}.
\ee

The goal of this work is to construct three dimensional solutions  with nonzero swirl  by perturbing background solutions.

Before we state the main result, some weighted H\"{o}lder norms are first introduced: For a bounded connected open set $\Om\subset \mbR^n$, let $\Ga$ be a closed portion of $\p\Om$. For ${\rm x,y} \in \Om$, set
\be\no
\de_{\rm x}\co \inf_{\rm z\in \Ga}|\rm x-\rm z|\q \text{and}\q \de_{\rm x,\rm y}\co \min(\de_{\rm x},\de_{\rm y}).
\ee
Given $k\in \mbR, \al\in (0,1)$, and $m\in \mb{Z}^+$, define the standard H\"{o}lder norms by
\be\no
\|u\|_{m,\Om}\co \sum_{0\leq |\beta|\leq m} \sup_{\rm x\in \Om} |D^{\beta}u(\rm x)|,\q [u]_{m,\al,\Om}\co \sum_{|\beta|=m} \sup_{{\rm x,y}\in\Om, x\neq y} \f{|D^{\beta}u(x)-D^{\beta}u(y)|}{|x-y|^{\al}}
\ee
where $D^{\beta}$ denotes $\p_{x_1}^{\beta_1}\cdots \p_{x_n}^{\beta_n}$ for a multi-index $\beta= (\beta_1,\cdots, \beta_n)$ with $\beta_j\in \mb{Z}_+$ and $|\beta|= \sum_{j=1}^n \beta_j$.
And, define weighted H\"{o}lder norms by
\be\no
\|u\|_{m,0,\Om}^{(k,\Ga)} &\co& \sum_{0\leq |\beta|\leq m} \sup_{\rm x\in\Om} \de_{\rm x}^{\max(|\beta|+k, 0)} |D^{\beta} u({\rm x})|,\\\no
[u]_{m,\al,\Om}^{(k,\Ga)} &\co& \sum_{|\beta|=m} \sup_{{\rm x,y}\in\Om, x\neq y}\de_{{\rm x,y}}^{\max(m+\al+k, 0)} \f{|D^{\beta}u({\rm x})-D^{\beta}u({\rm y})|}{|{\rm x-y}|^{\al}},\\\no
\|u\|_{m,\al,\Om} &\co& \|u\|_{m,\Om} + [u]_{m,\al,\Om},\q \|u\|_{m,\al,\Om}^{(k,\Ga)}\co \|u\|_{m,0,\Om}^{(k,\Ga)} + [u]_{m,\al,\Om}^{(k,\Ga)},
\ee
$C_{(k,\Ga)}^{m,\al}(\Om)$ denotes the completion of the set of all smooth functions whose $\|\cdot\|_{m,\al,\Om}^{(k,\Ga)}$ norms are finite.

\begin{problem}[Main problem]\lab{problem}
Given functions $(b, S_{en}, \mc{B}_{en}, \nu_{en}, \Phi_{bd}, p_{ex})$, find a solution $(\rho, {\bf u}, S, \Phi)\in [C^0(\ol{\mc{N}})\cap C^1(\mc{N})]^5\times { [C^1(\ol{\mc{N}})\cap C^2(\mc{N})]}$ to the nonlinear system \eqref{b17} with the following  boundary conditions
\be\lab{b131}
&&(S, \mc{B}, \La)(0, x_2,x_3)= (S_{en}, \mc{B}_{en}, r\nu_{en})(r)\qquad \text{on}\,\, \Ga_0,\\\lab{b131'}
&&{\bf u}\cdot {\bf e}_{r}=0 \qquad \text{on}\,\, \Ga_0,\\\lab{b132}
&&\Phi(x,r)= \Phi_{bd}(r)\qquad \text{on}\,\, \Ga_0\cup \Ga_L,\\\lab{b133}
&&{\bf u}\cdot {\bf n}_w = \p_{{\bf n}_w} \Phi=0 \qquad \text{on}\,\, \Ga_w,\\\lab{b134}
&&p(L, x_2,x_3)=p_{ex}(r)\q \q \text{on}\q \Ga_L,
\ee
where ${\bf n}_w$ is the unit normal vector field on $\Ga_w$ with being oriented interior to $\mc{N}$.
\end{problem}
\begin{remark}
For simplicity, we prescribe the Dirichlet boundary condition for $\Phi$ on $\Ga_0\cup\Ga_L$ as in \eqref{b132}. A physical condition such as $\p_x\Phi=E_{en}$ on $\Ga_0$ can be also considered by a simple adjustment of the argument in this paper.
\end{remark}

\smallskip
\begin{remark}\label{remark-14}
If an axisymmetric vector field ${\bf V}= V_x(x,r){\bf e}_x+ V_r(x, r){\bf e}_r+ V_{\th}(x, r) {\bf e}_{\th}$ is $C^1$ in $\mathbb{R}^3$, then it must satisfy $V_r(x,0)= V_{\th}(x, 0)\equiv 0$. From this point of view, we prescribe a compatibility conditions for $\nu_{en}$ as follows:
\be\lab{compatibility-nu}
\nu_{en}(0)=0.
\ee
\end{remark}
\bt\lab{main}
{\it Let $\bar{U}=(\bar{\rho}, \bar{\bf{u}}, S_0, \Phi_0)$ be the background solution in $\mc{N}$ associated with the entrance data $(b_0, J_0, S_0, \rho_0, E_0)$. Assume that $\nu_{en}$ satisfies \eqref{compatibility-nu}, and that
$\Phi_{bd}$ satisfies the compatibility condition
\be\lab{compatibility}\begin{array}{ll}
\p_{{\bf n}_w} \Phi_{bd}=0,\q \text{on}\q (\Ga_0\cup \Ga_L)\cap \ol{\Ga_w}.
\end{array}\ee
Given functions $(b, S_{en}, \mc{B}_{en}, \nu_{en}, \Phi_{bd}, p_{ex})(r)$, set
\be\lab{b213}\begin{array}{ll}
\om_1(b) \co \|b-b_0\|_{\al,\mc{N}},\\
\om_2(S_{en}, \mc{B}_{en}, \nu_{en})\co \|(S_{en}, \mc{B}_{en})- (S_0, \mc{B}_0 )\|_{1,\al,\Ga_0}+\|v_{en}\|_{1,\al,\Ga_0},\\
\om_3(\Phi_{bd}, p_{ex})\co \|\Phi_{bd}- \Phi_0\|_{2,\al,\Ga_0\cup\Ga_L}^{(-1-\al,\p(\Ga_0\cup\Ga_L))} + \|p_{ex}- p_L\|_{1,\al,\Ga_L}^{(-\al,\p\Ga_L)},
\end{array}\ee
and set
\begin{equation*}
\si:=\om_1(b)+\om_2(S_{en}, \mc{B}_{en}, \nu_{en})+\om_3(\Phi_{bd}, p_{ex})
\end{equation*}
for $p_L=p(\bar{\rho}(L), S_0)$.
Then there exists a constant $\si_1>0$ depending only on $(\ga, b_0, J_0, S_0, \rho_0, E_0, L,\al)$ so that if
\be\lab{b212}
\sigma\leq \si_1,
\ee
then {\rm Problem \ref{problem}} has an axially symmetric solution $(\rho, {\bf u}, S, \Phi)$ satisfying
\be\lab{e1}\begin{array}{ll}
\|(\rho,{\bf u})-(\bar{\rho}, \bar{{\bf u}})\|_{1,\al,\mc{N}}^{(-\al,\Ga_w)} +{{\|S-S_0\|_{1,\alp,\mc{N}}}}+ \|\Phi-\Phi_0\|_{2,\al,\mc{N}}^{(-1-\al,\Ga_w)}\le C\si,
\end{array}\ee
for the constant $C$ depending only on $(\ga, b_0, J_0, S_0, \rho_0, E_0, L,\al)$.
\medskip

Regarding $(S_{en}, \mc{B}_{en}-\Phi_{bd}, \nu_{en})$ as functions of $r\in[0,1]$, let us set $\om_4(S_{en}, \mc{B}_{en},\nu_{en},\Phi_{bd})$ as
\be\lab{e3}
\om_4(S_{en}, \mc{B}_{en},\nu_{en},\Phi_{bd})\co \|(S_{en}, \mc{B}_{en}-\Phi_{bd}, \nu_{en})- (S_0,\mc{B}_0, 0)\|_{C^{2,\alp}([0,1])}.
\ee
Then, there exists a constant $\si_2>0$ depending only on $(\ga, b_0, J_0, S_0, \rho_0, E_0, L,\al)$ such that if
\be\lab{e2}
\si+\om_4(S_{en}, \mc{B}_{en},\nu_{en},\Phi_{bd})\leq \si_2,
\ee
then the axially symmetric solution $(\rho, {\bf u}, S, \Phi)$ with satisfying \eqref{e1} is unique.
}\et

\begin{remark}
\label{remark-compay-SK}
Since it is assumed in \eqref{b213} that the axisymmetric functions $(S_{en}, \mc{B}_{en}, \nu_{en}, \Phi_{bd}, p_{ex})(r)$ are $C^1$ in $\Ga_0$, the compatibility conditions
\be\no
\p_r(S_{en}, \mc{B}_{en}, \nu_{en}, \Phi_{bd}, p_{ex})(0)=0
\ee
are naturally imposed.
\end{remark}

Unless otherwise specified, we say that {\it a constant $C$ is chosen depending only on the data} if $C$ is chosen depending only on $(\ga, b_0, J_0, S_0, \rho_0, E_0, L)$.
\medskip

In extending the results from \cite{bdx14} to three dimensional cases, the main difficulty is how to find a plausible Helmholtz decomposition of the velocity field. Fortunately, in the axisymmetric setting, the decomposition ${\bf u}= \na \varphi+\textit{curl }{\bf V}({\rm x})$ with ${\bf V}({\rm x})= h(x,r) {\bf e}_r+ \psi(x,r) {\bf e}_{\th}$ works. With using this representation, \eqref{b15} is decomposed as a weakly coupled system of second order elliptic equations for $(\varphi, \Phi,\psi)$,  and transport equations for $(S,\mc{B}, \La)$ for $\La=r\der_x h$.  In this reformulation, two new difficulties arise. As we shall see in the next section, the equation for $\psi$ contains a singular coefficient which blows up to infinity at $r=0$.
Secondly, due to nonzero swirl, a careful analysis of streamlines is required near the axis $r=0$ in order to solve the three transport equations $(S,\mc{B}, \La)$.

The rest of the paper is organized as follows.
In \S \ref{reformulation}, we introduce a Helmholtz decomposition for axisymmetric velocity fields, then reformulate the Euler-Poisson equations into a quasilinear second order elliptic system for  $(\varphi, \Phi, \psi)$, and three transport equations for the hyperbolic quantities $({S}, \mc{B}, \La)$. In \S \ref{linearization}, the unique solvability of a boundary value problem with a linearized elliptic system is discussed. We also prove the unique existence of $C^1$ solutions to transport equations. Finally, we implement an iteration to prove Theorem \ref{main}, the main result of this work, in \S \ref{iteration}.

\section{Helmholtz decomposition for axisymmetric flow of nonzero swirl}\lab{reformulation}

Assume that $(\rho, u_x, u_r, u_{\th}, S, \Phi)\in (C^1(\mc{N}))^5\times C^2(\mc{N})$ is a  solution to (\ref{b17}) with satisfying $\rho>0$ and $u_x>0$ in $\mc{N}$. We define a {\it pseudo-Bernoulli's} function $\mc{K}$ by
\be\lab{pseudoB}
\mc{K}\co \mc{B}-\Phi.
\ee
Then one can directly check that (\ref{b17}) is equivalent to the following system:
\be\lab{b170}
\begin{cases}
\p_x (\rho u_x) + \p_r (\rho u_r) + \f{\rho u_r}{r}=0,\\
\rho (u_x\p_x + u_r \p_r ) u_r - \f{\rho u_{\th}^2}{r} + \p_r p= \rho \p_r \Phi,\\
\rho (u_x\p_x + u_r \p_r ) \La = 0,\\
\rho (u_x\p_x + u_r \p_r ) \mc{K} = 0,\\
\rho (u_x\p_x + u_r \p_r ) S=0,\\
\Delta \Phi = \rho -b.
\end{cases}
\ee

As in \cite{bdx14}, we introduce a new representation of the velocity field ${\bf u}$. For
\[
{\bf V}({\rm x})= h(x,r){\bf e}_r+ \psi(x,r) {\bf e}_{\theta},\qquad \var({\rm x})=\var(x,r),
\]
set
\be\lab{Hdecomposition}
{\bf u}({\rm x})=\nabla \varphi({\rm x})+ \text{curl } {\bf V}({\rm x})
\ee
Suppose that $h{\bf e}_r$, $\psi {\bf e}_{\theta}$ and $\varphi$ are $C^2$ in $\mc{N}$. Then a straightforward computation yields
\be\lab{b110}\begin{array}{ll}
{\bf u}({\rm x})
&=\b(\f{1}{r}\p_r(r \psi) + \p_x \var\b) {\bf e}_x + (-\p_x \psi+ \p_r \var) {\bf e}_r + \p_x h(x,r){\bf e}_{\th},
\end{array}\ee
from which we derive that
\be\lab{b111}
u_x=\f 1r\p_r (r\psi) + \p_x \var,\q u_r=-\p_x \psi+\p_r \var, \q u_{\th}=\frac{\La}{r}=\p_x h.
\ee
\begin{remark}
\label{remark-psi-eqn}

The representation \eqref{b110} is well defined in $\mc{N}$ up to the axis $r=0$ if $u_r=u_{\theta}=0$ on $r=0$. Hereafter, any continuous vector field  ${\bf W}=W_x{\bf e}_x+W_r{\bf e}_r+W_{\theta}{\bf e}_{\theta}$ represented in the cylindrical coordinates is considered to satisfy $W_r=W_{\theta}=0$ on the axis $\{r=0\}$.

\end{remark}
Hereafter we denote the velocity field as
\be\lab{def-q}
{\bf u}={\bf q}(r, \psi, D\psi, D\var, \La),
\ee
for $D=(\p_x, \p_r)$
where ${\bf q}=(q_x, q_r, q_{\th})$ is given by the righthand sides of \eqref{b111}.

The vorticity field ${\bm \omega}({\rm x})=\text{curl }{\bf u}({\rm x})=\omega_x {\bf e}_x+ \omega_r {\bf e}_r+\omega_{\theta} {\bf e}_{\theta}$ is given by
\be\no
\om_x(x,r)=\f{1}{r} \p_{r}(r u_{\theta}),\q \om_r(x,r)=-\p_{x} u_{\theta},\q \om_{\th}(x,r)=\p_x u_r-\p_r u_x.
\ee

For $r>0$, we have
\be\no
\text{div} {\bf V}({\rm x})=\f{1}{r}\p_r(r h), \qquad
{\Delta }  {\bf V}({\rm x})=\Delta(h{\bf e}_r)+\Delta(\psi{\bf e}_{\theta})=\b(\Delta-\f{1}{r^2}\b)h {\bf e}_r+ (\Delta-\f{1}{r^2})\psi {\bf e}_{\th}.
\ee

Substituting the representations into
$
\text{curl }{\bf u}({\rm x})=\text{curl }\text{curl } {\bf V}({\rm x})=\na \text{div} {\bf V}({\rm x})-{\Delta}  {\bf V}({\rm x}),
$
we obtain that
\be\label{vector-eq-t}
-\Delta(\psi {\bf e}_{\theta})=\om_{\theta}{\bf e}_{\theta}.
\ee

And, from \eqref{pseudoB}, (\ref{b170}) and \eqref{b111} , it follows that
\be\lab{b171}\begin{array}{ll}
\om_{\theta}&= \p_x u_r- \p_r u_x=\f{1}{u_x}[(u_x\p_x+u_r\p_r) u_r-\p_r\f{1}{2}|{\bf u}|^2+u_{\th}\p_r u_{\th}]\\
&=\f{1}{\f{1}{r}\p_r(r\psi)+ \p_x \var}\b(T(\rho, S)\p_r S-\p_r \mc{K}+\f{\La}{r^2}\p_r \La\b),
\end{array}\ee
for
\be\no
T(\rho, S)= \f{A}{c_v (\ga-1)} \textit{exp}\b(\f{S}{c_v}\b) \rho^{\ga-1}.
\ee
\begin{remark}\label{remark-psi-regularity}
Note that if $\psi{\bf e}_{\theta}$ is $C^2$ in $\mc{N}$, and if $\om_{\theta}=0$ on $r=0$, the equation \eqref{vector-eq-t} is well defined. Since $\psi$ is assumed to depend only on $(x,r)$ in the cylindrical coordinates, if $\lim_{r\to 0+}(\Delta-\frac{1}{r^2})\psi(x,r)=\lim_{r\to 0+} \left(\psi_{xx}+\frac 1r\p_r(r\psi_r)-\frac{1}{r^2}\psi\right)(x,r)=0$, then the equation $-\Delta(\psi {\bf e}_{\theta})=\om_{\theta}{\bf e}_{\theta}$ can be simply written as
\be\lab{b112}
-\b(\Delta - \f{1}{r^2}\b) \psi(x,r)=\om_{\theta}(x,r).
\ee
The main challenge of this work is to find a solution $\psi$ to \eqref{b112} with satisfying 
\be \label{limits-psi}
\lim_{r\to 0+}(\Delta-\frac{1}{r^2})\psi(x,r)=\lim_{r\to 0+} \left(\psi_{xx}+\frac 1r\p_r(r\psi_r)-\frac{1}{r^2}\psi\right)(x,r)=0\quad\tx{for all}\quad x\in[0,L]
\ee
provided that the compatibility condition stated in Remark \ref{remark-compay-SK} hold. If the existence of such a solution $\psi$ is proved, then $\psi(x,r){\bf e}_{\theta}$ solves the vector equation \eqref{vector-eq-t}. Details are discussed in Proposition \ref{bl510}.
\end{remark}

It follows from (\ref{b16}) and \eqref{pseudoB} that
\be\lab{b139}
\rho=H(S, \mc{K}+\Phi-\f 12 |{\bf q}(r, \psi, D\psi, D\var, \La)|^2),
\ee
where
\be
\label{definition-H}
H(\xi,\tau)=\b[\f{(\ga-1)\tau}{\ga A}\textit{exp }\b(-\f{\xi}{c_v}\b)\b]^{\f{1}{\ga-1}}.
\ee

Using \eqref{b171} and \eqref{b112}, it can be directly checked that (\ref{b170}) is equivalent to the following system:
\be\lab{b141}
&&\text{div }\b(H(S, \mc{K}+\Phi-\f 12 |{\bf q}(r, \psi, D\psi, D\var, \La)|^2){\bf q}(r, \psi, D\psi, D\var, \La)\b)=0,\\\lab{b142}
&&\Delta \Phi= H(S, \mc{K}+\Phi-\f 12 |{\bf q}(r, \psi, D\psi, D\var, \La)|^2)-b,\\\lab{b143}
&&- \Delta (\psi{\bf e}_{\theta}) = \f{T(\rho, S) \p_r S- \p_r \mc{K}+\f{\La}{r^2}\p_r \La}{\f 1r \p_r(r \psi) + \p_x \varphi}{\bf e}_{\theta},\\\lab{b144}
&&\rho \b[\b(\f 1r \p_r (r\psi)+ \p_x \varphi\b)\p_x + (-\p_x \psi + \p_r \varphi) \p_r \b] S=0,\\\lab{b145}
&&\rho \b[\b(\f 1r \p_r (r\psi)+ \p_x \varphi\b)\p_x + (-\p_x \psi + \p_r \varphi) \p_r \b] \mc{K}=0,\\\lab{b146}
&&\rho \b[\b(\f 1r \p_r (r\psi)+ \p_x \varphi\b)\p_x + (-\p_x \psi + \p_r \varphi) \p_r \b] \La=0.
\ee
\medskip

Next, we derive the corresponding boundary conditions for $(\var,\psi)$ from \eqref{b131}, \eqref{b133} and \eqref{b134}.
If $\var=k_1$ for some constant $k_1$ and $\psi_x=0$ on $\Ga_0$, then ${\bf u}$ given by \eqref{b110} satisfies \eqref{b131} on $\Ga_0$. Also, if $\psi=k_2$ for some constant $k_2$ and $\der _{{\bf n}_w}\var=0$  on $\Ga_w$, then \eqref{b133} holds. So we prescribe:
\be\lab{b135}
\varphi=0,\q \q \text{on}\q \Ga_0, \q  \q \p_{{\bf n}_w} \varphi=0\q\q \text{on}\q \Ga_{w},
\ee
\be\lab{b137}
-\p_x \psi(0,r)=0\q\q \text{on}\q \Ga_0\q \q \text{and}\q  \psi=0\q\q \text{on}\q \Ga_w.
\ee
At the exit $\Ga_L$ of nozzle $\mc{N}$, we fix a boundary condition for $\psi$ as
\be\lab{b138}
\p_x \psi(L,r)=0\q \q \text{on}\q \Ga_L.
\ee
We also require
\be\lab{b138'}
\psi(x,0)=0,\q \forall x\in [0, L]
\ee
so that ${\bf u}(\rm x)$ given by \eqref{b110}  satisfies a necessary condition to be a $C^1$ axisymmetric vector field in $\mc{N}$. See Remark \ref{remark-14}.

We collect the boundary conditions for $(\var, \psi,\Phi, S, \mc{K}, \La)$ as follows
\be\lab{b317}
&&\varphi=0\q \text{on}\q \Ga_0\q \text{and}\q \p_{{\bf n}_w}\var =0\q \text{on}\q \Ga_w,\\\lab{b318}
&&\Phi= \Phi_{bd} \q \text{on}\q \Ga_0\cup\Ga_L,\q \p_{{\bf n}_w} \Phi=0\q \text{on}\q \Ga_w,\\\lab{b319}
&&\p_{x}\psi(0,r)=0\q \text{on}\q \Ga_0,\q \p_x \psi(L,r)=0\q \text{on}\q \Ga_L, \q \text{and}\q \psi=0\q \text{on}\q \Ga_w,\\\lab{b320}
&&(S, \mc{K}, \La)=(S_{en}, \mc{K}_{en}, r\nu_{en})\q \text{on}\q \Ga_0,\\\lab{b320'}
&& A \text{exp}\b(\f{S}{c_v}\b) H^{\ga}\b(S, \mc{K}+\Phi-\f 12 |{\bf q}(r, \psi, D\psi, D\var, \La)|^2\b)= p_{ex}\q \text{on}\q \Ga_L.
\ee

Theorem \ref{main} will directly follow from the following theorem.
\bt\lab{main1}
{\it  Let $\sigma>0$ be defined by \eqref{b213} in Theorem \ref{main}.
Under the same assumptions as {\rm Theorem \ref{main} }, there exists a constant $\si_3>0$ depending only on the data and $\al$ so that if
\be\lab{b321}
\si\leq \si_3,
\ee
then the boundary value problem (\ref{b141})-(\ref{b146}) with (\ref{b138'})-(\ref{b320'}) has an axially symmetric solution $(\var, \psi,\Phi, S, \mc{K}, \La)$ satisfying
\be\lab{b322-1}\begin{array}{ll}
{{\|(\var, \Phi)-(\var_0, \Phi_0)\|_{2,\al,\mc{N}}^{(-1-\al,\Ga_w)} +\|\psi {{\bf e}}_{\theta}\|_{2,\alp,\mc{N}}+ \|(S,\mc{K},\La)-(S_0, \mc{B}_0, 0)\|_{1,\al,\mc{N}}}}\le C\sigma,
\end{array}\ee
for constant $C>0$ depending only on the data and $\al$.

Moreover, there exists a constant $\si_4>0$ depending only on the data and $\al$ such that if
\be\lab{b323}
\si+\om_4(S_{en}, \mc{B}_{en},\nu_{en},\Phi_{bd})\leq \si_4,
\ee
then the axially symmetric solution $(\var, \psi,\Phi, S, \mc{K}, \La)$ with satisfying \eqref{b322-1} is unique.
}\et

\section{Linear boundary value problem associated with  (\ref{b141})-(\ref{b146})}
\lab{linearization}

\subsection{Linearization of (\ref{b141})-(\ref{b146})}
To prove Theorem \ref{main1}, we first investigate the unique solvability of a linear boundary value problem associated with the nonlinear boundary value problem of \eqref{b141}--\eqref{b146} and \eqref{b138'}--\eqref{b320'}.
For $(\zeta,\eta, z)\in\mbR^3$, ${\bf v}=({\rm v}_1,{\rm v}_2,{\rm v}_3), {\bf s}=(s_1,s_2,s_3)\in \mbR^3$, define ${\bf A}=(A_1, A_2, A_3)$ and $B$ by
\be\no\begin{array}{ll}
A_j(\zeta, \eta, z, {\bf v}, {\bf s}) = B(\zeta,\eta,z, {\bf v}, {\bf s}) {\rm v}_j\q\q\q \textit{for}\q j=1,2,3,\\
B(\zeta,\eta,z, {\bf v}, {\bf s}) =H\b(\zeta,\eta+z-\f 12|{\bf v}+{\bf s}|^2\b).
\end{array}\ee
Set $\mc{K}_0=\mc{B}_0$ for $\mc{B}_0$ given by \eqref{def-B0}, and denote ${\bf V}_0=(S_0, \mc{K}_0, \Phi_0, D \var_0, {\bf 0})$. Set
\be\lab{b44}
 a_{ij}(x)= \p_{{\rm v}_j} A_i({\bf V}_0),\q b_i(x)= \p_z A_i({\bf V}_0)\q c_i(x)=\p_{{\rm v}_i} B({\bf V}_0),\q d(x)= \p_z B({\bf V}_0)
\ee
for $i,j=1,2,3$. Note that $a_{ij}, b_i, c_i, d$ are functions of $x$ only because ${\bf V}_0$ depends only on $x\in[0,L]$.

\bl\lab{bl45}
{\it The coefficients $a_{ij}, b_i, c_i, d$ defined by (\ref{b44}) satisfy the following properties:
\begin{enumerate}[(a)]
  \item The matrix $[a_{ij}(x)]_{i,j=1}^3$ is diagonal and $a_{22}(x)\equiv a_{33}(x)$ in $\mc{N}$, and there exists a constant $\nu_1>0$ satisfying
  \be\lab{b45}
  \nu_1 I_3\leq [a_{ij}(x)]_{i,j=1}^3 \leq \f{1}{\nu_1} I_3\q \textit{for all}\q {\bf x}\in \mc{N},
  \ee
  where the constant $\nu_1$ depends only on the data.
  \item For each $i=1,2, 3$, we have
  \be\lab{b47}
  b_1(x)+ c_1(x)\equiv 0,\q\q b_j(x)=c_j(x)\equiv 0\q \textit{for}\,\,j=2,3\q\textit{in}\q\mc{N}.
  \ee
  \item For each $k\in \mb{Z}_+$, there exists a constant $C_k>0$ depending only on the data and $k$ such that
  \be\lab{b46}
  \sum_{i=1}^3 \|a_{ii}\|_{k,\mc{N}}+ (\|b_1\|_{k,\mc{N}}+\|c_1\|_{k,\mc{N}})+ \|d\|_{k,\mc{N}}\leq C_k.
  \ee
  \item There exists a constant $\nu_2>0$ depending only on the data such that
  \be\lab{b48}
  d(x)\geq \nu_2\q \textit{in}\q\mc{N}.
  \ee
\end{enumerate}
}\el
Lemma \ref{b145} has been proved in \cite{bdx14} so we skip to prove it. But we remark that \eqref{b47} in Lemma \ref{b145} is essential for obtaining the well-posedness of the linearized elliptic system (\ref{b24}) stated below.

For ${\bf q}(r, \psi, D\psi, D\var, \Lambda)$ defined by \eqref{def-q}, set
\[
{\bf t}(r, \psi,D\psi, \La):={\bf q}(r, \psi, D\psi, D\var, \La)-\na\var.
\]
For later use, we further represent ${\bf t}(r, \psi,D\psi, \La)$ as
\begin{equation*}
{\bf t}(r, \psi,D\psi, \La)={\bf t}_1(r,\psi, D\psi)+{\bf t}_2(r, \La)
\end{equation*}
for
\begin{equation}
\label{definition-tvecs}
{\bf t}_1(r,\psi, D\psi)=\frac 1r\p_r(r\psi){\bf e}_x-\p_x\psi{\bf e}_r,\quad\quad
{\bf t}_2(r, \La)=\frac{\La}{r} {\bf e}_{\theta}.
\end{equation}

Then (\ref{b141}) and (\ref{b142}) can be written as
\be\no\begin{array}{ll}
\text{div}\left({\bf A}(S,\mc{K},\Phi,D\var, {\bf t}(r, \psi,D\psi, \La))\right)=-
\text{div}\left(B(S,\mc{K},\Phi,D\var, {\bf t}(r, \psi, D\psi, \La))){\bf t}(r, \psi, D\psi, \La)\right),\\
\Delta \Phi=B(S,\mc{K},\Phi,D \var, {\bf t}(r, \psi, D\psi, \La))-b.
\end{array}\ee
For $(\Phi_0, \var_0)$ given by \eqref{background5}, set
\be\no
(\Psi, \phi)\co  (\Phi- \Phi_0, \var-\var_0),
\ee
Set $\mc{K}_0:=\mc{B}_0$, and denote
\be\no
\mc{U}\co (\Psi, \phi,\psi),\q \mc{W}\co (S,\mc{K},\La),\q \mc{W}_0\co (S_0,\mc{K}_0, 0),\q \mc{W}_{en}\co (S_{en}, \mc{K}_{en}, r\nu_{en}).
\ee
Then $(\Psi,\phi)$ satisfy the equations
\be\lab{b23}\begin{cases}
L_1(\Psi,\phi)= \text{div }{\bf F}({\bf x}, \mc{W}-\mc{W}_0, \Psi, D \phi,{\bf t}(r, \psi, D\psi, \La)),\\
L_2(\Psi,\phi)= f_1({\bf x}, \mc{W}-\mc{W}_0, \Psi, D \phi, {\bf t}(r, \psi, D\psi, \La)),
\end{cases}\ee
where $L_1, L_2$, ${\bf F}=(F_1, F_2, F_3)$ and $f_1$ are defined as follows:
\be\lab{b24}\begin{cases}
L_1(\Psi, \phi)=  \sum_{i=1}^3 \p_i(a_{ii}(x)\p_i \phi)+ \p_x (b_1(x)\Psi),\\
L_2(\Psi, \phi)=\Delta \Psi-c_1(x)\p_x \phi- d(x) \Psi,
\end{cases}\ee
and
\be\lab{b25}
F_i({\bf x},Q) &=&- B({\bf V}_0+ Q){\bf s}_i -\int_0^1 D_{(\eta_1,\eta_2,{\bf s})} A_i({\bf V}_0+t Q) dt\cdot (\eta_1,\eta_2, {\bf s})\\\no
&\q& -\int_0^1 D_{(z,{\bf v})} A_i({\bf V}_0+ t Q)- D_{(z,{\bf q})} A_i({\bf V}_0) dt\cdot (z,{\bf v}),\q i=1,2,3,\\\lab{b26}
f_1({\bf x},Q) &=& \int_0^1 D_{(\eta_1,\eta_2,{\bf s})} B({\bf V}_0+ t Q) dt\cdot (\eta_1,\eta_2, {\bf s})- (b-b_0)\\\no
&\q& + \int_0^1 D_{(z,{\bf v})} B({\bf V}_0+t Q)- D_{(z,{\bf v})} B({\bf V}_0) dt \cdot (z, {\bf v})
\ee
with $Q= (\eta_1, \eta_2, z, {\bf v}, {\bf s})\in \mbR^3\times (\mbR^3)^2$.
We subtract the expression
\be\no
A \textit{exp}\b(\f{S_0}{c_v}\b) H^{\ga}\b(S_0,\mc{K}_0+\Phi_0-\f12 |D \var_0|^2\b)=p_L \quad \textit{on}\q \Ga_L
\ee
from (\ref{b134}) to get
\be\lab{b27}
\p_x \phi(L,r)= g(r,\mc{W}-\mc{W}_0, D \phi, {\bf t}(r, \psi, D\psi, \La))\q \textit{on}\q \Ga_L
\ee
with $g$ defined by
\be\no
g(r, \mc{K}-\mc{K}_0, D \phi, {\bf t}(r, \psi, D\psi, \La)) &=& -\f{1}{r}\p_r(r\psi)+ \f{\mc{K}-\mc{K}_{0}+\Psi_{bd}-\f 12|D \phi+{\bf t}(r, \psi, D \psi, \La)|^2}{\bar{u}(L)}\\\lab{b28}
&\q&-\f{p_{ex}^{\f{\ga-1}{\ga}}\b(A\textit{exp}\b(\f{S}{c_v}\b)\b)^{\f{1}{\ga}}-p_L^{\f{\ga-1}{\ga}}\b(A\textit{exp}\b(\f{S_0}{c_v}\b)\b)^{\f1{\ga}}}{(\ga-1)\bar{u}(L)}
\ee
for $\Psi_{bd}= \Phi_{bd}-\Phi_0$.
From (\ref{b132}) and (\ref{b135}), we can derive the boundary conditions for $\Psi$ and $\phi$:
\be\lab{b212}
&&\phi=0\q \textit{on}\q \Ga_0,\q\q \p_{{\bf n}_w} \phi=0\q \textit{on}\q \Ga_w,\\\lab{b213-new}
&&\Psi=\Psi_{bd}\q \textit{on}\q \Ga_0\cup\Ga_L,\q \q \p_{{\bf n}_w}\Psi=0\q \textit{on}\q \Ga_w.
\ee
Since $\mc{W}_0$ is a constant vector, (\ref{b143}) can be rewritten as
\be\lab{b210}
-\Delta(\psi{\bf e}_{\theta})= f_2\left({\rm x}, \mc{W}-\mc{W}_0,\Psi, D \phi,{\bf t}(r, \psi,D\psi, \La), \p_r(\mc{W}-\mc{W}_0)\right) {\bf e}_{\theta}
\ee
for $f_2$ defined by
\be\lab{b211}
f_2({\rm x},Q, \p_r(\mc{W}-\mc{W}_{0}))= \f{T(B({\bf V}_0+Q), S_0+\eta_1)\p_r(S-S_0)-\p_r(\mc{K}-\mc{K}_0)+\f{\Lambda}{r^2}\p_r(\La-\La_0)}{\p_x \var_0(x)+ q_1 +s_1}.
\ee
By introducing the vector field
\be\lab{b251}
\begin{split}
&{\bf M}(S, \mc{K}, \Psi, \nabla\phi, {\bf t}(r, \psi,D\psi, \La))\\
&:=H(S, \mc{K}+\Phi_0+ \Psi-\f 12 |\nabla \varphi_0+\na \phi+{\bf t}(r, \psi,D\psi, \La)|^2)\left(\nabla \varphi_0+\na \phi+{\bf t}(r, \psi,D\psi, \La)\right)\\
&=H(S, \mc{K}+\Phi_0+ \Psi-\f 12 |\nabla \varphi+{\bf t}(r, \psi, D\psi, \La)|^2)\b[\b(\p_x\var+ \f 1r\p_r(r\psi)\b){\bf e}_{x}+(\p_r\var-\p_x \psi){\bf e}_r\b],
\end{split}
\ee
we can rewrite the transport equations (\ref{b144})--(\ref{b146}) in the form of
\be\lab{b252}\begin{array}{ll}
{\bf M}\cdot \na \mc{W}= 0,\q &\text{in}\q \mc{N},\\
\mc{W}(0, r)=\mc{W}_{en}(r),\q &\text{on}\q \Ga_0.
\end{array}\ee
Here ${\bf M}=M_x{\bf e}_x + M_r {\bf e}_r$ satisfies
\be\lab{b253}
\p_x(r M_x)+ \p_r (r M_r)=0\q \q \forall (x,r)\in \Om=[0,L]\times [0,1],\q M_r(x,0)=M_r(x,1)=0.
\ee
Therefore \eqref{b252} can be regarded as transport equations in a two dimensional rectangular domain with the divergence-free vector field $r{\bf M}$.
\medskip

\subsection{Linearized elliptic system for $(\Psi, \phi)$}

Suppose that $\mc{F}=(\mc{F}_1,\mc{F}_2, \mc{F}_3)(x,r)\in [C_{(-\al,\Ga_w)}^{1,\al}(\mc{N})]^3$, $\mathbbm{f}_1(x,r)\in C^{\al}(\mc{N})$ and $\mathbbm{g}(r)\in C^{1,\al}_{(-\al,\p\Ga_L)}(\Ga_L)$. Consider the linear system
\be\lab{b41}\begin{cases}
L_1(\Psi,\phi)= \textit{div}\mc{F}\\
L_2(\Psi,\phi)= \mathbbm{f}_1
\end{cases}\q \textit{in}\q \mc{N}\ee
with boundary conditions
\be\lab{b42}
\phi=0\q \textit{on}\q \Ga_0,\q \p_{{\bf n}_w} \phi=0\q \textit{on}\q \Ga_w,\q \p_x \phi=\mathbbm{g}(r)\q \textit{on}\q \Ga_L\\\lab{b43}
\Psi=\Psi_{bd}\q\textit{on}\q \Ga_0\cup\Ga_L,\q \q \p_{{\bf n}_w} \Psi=0\q \textit{on}\q \Ga_w.
\ee

\bl\lab{bl42}
{\it Suppose that $\mc{F}= (\mc{F}_1, \mc{F}_2, \mc{F}_3)(x,r)\in (C^{1,\al}_{(-\al,\Ga_w)}(\mc{N}))^3$, $\mathbbm{f}_1(x,r)\in C^{\al}(\ol{\mc{N}})$ and $\mathbbm{g}(r)\in C^{1,\al}_{(-\al,\p\Ga_L)}(\Ga_L)$ for $\al\in (0,1)$. If, in addition, $\Psi_{bd}$ satisfies the compatibility condition
\be\lab{b410}
\p_{{\bf n}_w} \Psi_{bd}=0\q \text{on} \q(\ol{\Ga_0}\cup\ol{\Ga_L}) \cap \ol{\Ga_w},
\ee
then the linear boundary value problem (\ref{b41})-(\ref{b43}) has a unique axially symmetric solution $(\phi, \Psi)=(\phi, \Psi)(x,r)\in (C^{1,\al}(\ol{\mc{N}})\cap C^{2,\al}(\mc{N}))^2$. Moreover, $(\phi, \Psi)$ satisfy the estimate
\be\lab{b411}
\begin{split}
&
\|(\phi,\Psi)\|_{1,\al,\mc{N}} \leq C_1(\|\mathbbm{g}\|_{\alp,\Ga_L}+\|\Psi_{bd}\|_{1,\al,\mc{N}}+\|\mc{F}\|_{\al,\mc{N}}+\|\mathbbm{f}_1\|_{\al,\mc{N}}),
\\
&\|(\phi,\Psi)\|_{2,\al,\mc{N}}^{(-1-\al,\Ga_w)} \leq C_1(\|\mathbbm{g}\|_{1,\al,\Ga_L}^{(-\al,\p\Ga_L)}+\|\Psi_{bd}\|_{2,\al,\mc{N}}^{(-1-\al,\Ga_w)}+\|\mc{F}\|_{1,\al,\mc{N}}^{(-\al,\Ga_w)}+\|\mathbbm{f}_1\|_{\al,\mc{N}})
\end{split}
\ee
where $C_1>0$ depends only on the data and $\al$.
}\el

\bpf
The well-posedness of \eqref{b41}--\eqref{b43} and the estimate \eqref{b411} have been proved in \cite[Proposition 4.1]{bdx}. So we only prove the axi-symmetric property of the unique solution $(\phi,\Psi)$ to \eqref{b41}--\eqref{b43}.

 For any $\theta\in [0,2\pi)$, define
\be\no
f^{\theta}({\bf x})= f(x_1, x_2 \cos\theta- x_3 \sin\theta, x_2 \sin \theta+x_3\cos\theta).
\ee
Since $a_{22}\equiv a_{33}$ and $\p_{x_2}^2+\p_{x_3}^2$ is invariant under the rotation group, we have
\be\no\begin{array}{ll}
L_1(\Psi^{\theta}, \phi^{\theta})= \text{div}\mc{F}(x_1, x_2 \cos\theta- x_3 \sin\theta, x_2 \sin \theta+x_3\cos\theta)=\text{div}\mc{F}({\rm x}),\\
L_2(\Psi^{\theta}, \phi^{\theta})= \mathbbm{f}_1(x_1,x_2 \cos\theta- x_3 \sin\theta, x_2 \sin \theta+x_3\cos\theta)= \mathbbm{f}_1({\rm x}),
\end{array}\ee
where we have used the axially symmetric properties of $\mc{F}$ and $\mathbbm{f}_1$. Therefore $(\phi^{\theta}, \Psi^{\theta})$ is also a solution to \eqref{b41}--\eqref{b43}. By the uniqueness of a solution, we conclude that $(\phi^{\theta}, \Psi^{\theta})=(\phi, \Psi)$, therefore
 $\phi$ and $\Psi$ are axially symmetric.
\epf

\subsection{Elliptic equation for $\psi$ with a singular coefficient}

We consider the following boundary value problem for a vector field ${\bf V}: \ol{\mc{N}}\rightarrow \R^3$:
\be\lab{51-vec}\begin{array}{ll}
-\Delta\,{\bf V}= \mathbbm{f}_2(x,r){\bf e}_{\theta}\q &\text{in}\q \mc{N},\\
\p_x {\bf V}= {\bf 0}\q &\text{on}\q \Ga_0\cup \Ga_L,\\
\bf V=0\q &\text{on}\q \Ga_w.
\end{array}\ee
If $\mathbbm{f}_2{\bf e}_{\theta}$ is $C^{\alp}$ in $\mc{N}$, then the standard elliptic theory(\cite{gt}) yields that \eqref{51-vec} has a unique solution ${\bf V}: \ol{\mc{N}}\rightarrow \R^3$ which satisfies the estimate
\[
\|{\bf V}\|_{2,\alp,\mc{N}}\le C\|\mathbbm{f}_2{\bf e}_{\theta}\|_{\alp, \mc{N}}
\]
for a constant $C>0$ depending only on $L$ and $\alp$. Note that the continuity of $\mathbbm{f}_2(x,r){\bf e}_{\theta}$ in $\mc{N}$ naturally implies that the function $\mathbbm{f}_2(x,r)$ satisfies the compatibility condition
\be\lab{b511}
\mathbbm{f}_2(x,0)\equiv 0,\q \forall x\in [0, L].
\ee

As discussed in Remark \ref{remark-psi-eqn}, we will show that the unique solution ${\bf V}$ to \eqref{51-vec} has the form of 
\begin{equation}
\label{estimate-V}
{\bf V}=\psi(x,r){\bf e}_{\theta},
\end{equation}
where $\psi$ solves
\be\lab{51}\begin{array}{ll}
-\b(\p_{xx}+\frac 1r\p_r(r\p_r) -\f1{r^2}\b) \,\psi= \mathbbm{f}_2(x,r)\q &\text{in}\q \mc{N},\\
-\p_x \psi(0,r)= 0\q &\text{on}\q \Ga_0,\\
\psi=0\q &\text{on}\q \Ga_w,\\
\p_x\psi=0\q &\text{on}\q \Ga_L,\\
\psi=0\q &\text{on}\q \mc{N}\cap \{r=0\}
\end{array}\ee
in the following sense:
\begin{itemize}
\item[(i)] As a function of $(x,r)$ in a two dimensional rectangle $\Om=(0,L)\times (0,1)$, $\psi$ is $C^2$ in $\Om$, and satisfies the equation and all the boundary conditions of \eqref{51} pointwisely;

\item[(ii)] As a function of ${\bf x}\in \mc{N}$, $\psi$ is not necessarily $C^2$ up to $r=0$, but it is a solution to \eqref{51} in distribution sense, where we write as $(\Delta-\frac{1}{r^2})\psi=\b(\p_{xx}+\frac 1r\p_r(r\p_r) -\f1{r^2}\b) \,\psi$;

\item[(iii)] As a vector field in $\mc{N}$, $\psi{\bf e}_{\theta}$ is a classical solution to \eqref{51-vec}.
\end{itemize}

\bp\lab{bl510}
{\it
For a fixed $\al\in(0,1)$, suppose that a vector field $\mathbbm{f}_2(x,r){\bf e}_{\theta}:\ol{\mc{N}}\rightarrow \R^3$ is $C^{\al}$ in $\mc{N}$. Note that the compatibility condition \eqref{b511} is naturally imposed.  Then the linear boundary value problem \eqref{51-vec} has a unique solution ${\bf V}:\ol{\mc{N}}\rightarrow \R^3$ which satisfies the estimate \eqref{estimate-V}. Furthermore, ${\bf V}$ is represented as
\begin{equation}
\label{V-representation}
\bf{V}(\rm x)=\psi(x,r){\bf e}_{\theta}\quad\tx{in}\quad\ol{\mc{N}},
\end{equation}
where $\psi$ solves \eqref{51} in the sense of (i)--(iii) stated above.
Furthermore, $\psi$ satisfies the estimate
\be\lab{b512}
\|\psi\|_{2,\al,\Om} \leq C_2\|\mathbbm{f}_2\|_{\al,\mc{N}}
\ee
for a constant $C_2>0$ depending only on the data and $\al$, where $\psi$ is regarded as a function defined in the two dimensional rectangle $\Om=[0, L]\times [0,1]$.
Note that ${\bf V}$ is well defined by \eqref{V-representation} due to the condition $\psi(x,0)=0$ for all $x\in[0, L]$.

}
\ep

\bpf
We prove this proposition in two methods.
\medskip

{\emph{(Method I)}}
{\textbf{1.}} In order to represent ${\bf V}$ in the form of \eqref{V-representation}, we need to find a solution $\psi$ to \eqref{51}.

The main idea to solve \eqref{51} is to rewrite it as a boundary value problem in $\mathbb{R}^5$ so that the singular term $\frac{\psi}{r^2}$ is removed from the equation for $\psi$. This idea has been used extensively in the study of Navier-Stokes equations, see {\cite{knss,lw09}}.
\medskip

Set
\be\label{def-xi-psi}
\xi(x,r):=\f{\psi(x,r)}{r},\q f(x,r):=-\f{\mathbbm{f}_2(x,r)}{r}.
\ee
We regards $\xi$ and $f$ as functions defined in
\be\no
\mc{D}\co(0, L)\times \{y'\in \mbR^4: |y'|< 1\}\subset \mathbb{R}^5,
\ee
where ${\bf y}=(x,r, {\bm \omega})\in \mathbb{R}\times \mathbb{R}^+\times S^3$ represent cylindrical coordinates in $\mathbb{R}^5$.
By the compatibility condition \eqref{b511}, we have
\begin{equation}
\label{f-estimate}
|f(x,r)|\le [\mathbbm{f}_2]_{\al} r^{-1+\al}.
\end{equation}
Define
\be\no
{\bf F}(x, y')=(0, F(x, r)y_2, F(x,r)y_3, F(x,r) y_4, F(x,r) y_5),\q\forall (x,y')\in\mc{D},
\ee
with
\be\no
F(x,r)=\f{1}{r^4} \int_0^r s^3 f(x,s) ds,
\ee
so that $f(x, r)= \text{div}_{{\bf y}} {\bf F}(x, y')$ for $\forall (x,y')\in \mc{D}$. By using \eqref{f-estimate}, one can directly check that
${\bf F}\in C^{\al}(\ol{\mc{D}})$ and
\be\lab{force}
\|{\bf F}\|_{\al,\ol{\mc{D}}}\leq C\|\mathbbm{f}_2\|_{\al,\ol{\mc{N}}}.
\ee
A formal computation shows that $\psi$ solves \eqref{51} if  $\xi$ solves
\be\lab{b531}\begin{array}{ll}
\Delta_{{\bf y}} \xi= \text{div}_{{\bf y}} {\bf F}(x, y') \q &\text{in}\q \mc{D},\\
-\p_x \xi(0,y')= 0\q &\text{on}\q B_0\co\{0\}\times\{y'\in \mbR^4: |y'|\leq 1\} ,\\
\xi(x, y')=0\q &\text{on}\q B_w\co [0,L]\times \{y': |y'|=1\},\\
\p_x\xi(L, y')=0\q &\text{on}\q B_{L}\co \{L\}\times\{y'\in \mbR^4: |y'|\leq 1\}.
\end{array}\ee
\eqref{b531} has a unique weak solution $\xi\in H^1(\mc{D})$, and the weak solution satisfies
\be\lab{536}
\|\xi\|_{1,\al,{\mc{D}}}\leq C\|{\bf F}\|_{\al,{\mc{D}}}\leq C\|\mathbbm{f}_2\|_{\al, {\mc{N}}}.
\ee
The estimate \eqref{536} is obtained by adjusting Theorem 3.13 of \cite{hl}. 
{As in Lemma \ref{bl42}, we can prove that $\xi$ is axially symmetric (i.e. $\xi({\bf y})=\xi(x, |y'|)$) by using the special orthogonal group $SO_4$ and the uniqueness of a weak solution to (\ref{b531}).}
Due to the uniqueness of a weak solution to \eqref{b531}, $\xi$ is axially symmetric i.e. $\xi({\bf y})=\xi(x, r)$.

For each constant $\de\in(0,1)$, define
$\mc{D}_{\de}=\{{\bf y}\in \mc{D}: r > \de\}$. Since $f$ is $C^{\al}$ in $\mc{D}$ away from $r=0$, the standard Schauder estimate(\cite{gt}) yields a constant  $C_{\delta}>0$ depending on $(\delta, \al)$ to satisfy

\be\lab{b535}\begin{array}{ll}
\|\xi\|_{2,\al,\mc{D}_{\de}}&\leq C_{\de} \|f\|_{\al, \mc{D}_{\de/2}}\leq C_{\de} \|\mathbbm{f}_2\|_{\al, \mc{N}}.
\end{array}\ee

{Therefore $\psi=r\xi$ satisfies all the boundary conditions in \eqref{51} and the equation $-(\Delta-\frac{1}{r^2})\psi=\mathbbm{f}_2(x,r)$ in $\mc{N}\setminus\{r=0\}$ in the classical sense.}
\smallskip

\textbf{2.} Next, we show that $\psi$ is $C^{2}$ with respect to $(x,r)$ especially up to $r=0$. Regarding $\psi$ as a function of two variables $(x,r)\in \Om$ for $\Om=[0,L]\times [0,1]$, $\psi$ solves the following two dimensional linear boundary value problem:
\be\lab{540}\begin{array}{ll}
(\p_x^2+\p_r^2) \psi=-\p_r \xi- \mathbbm{f}_2\left(= 
-\b(\f 1 r\p_r-\f1{r^2}\b)\psi-\mathbbm{f}_2\right), \q &\text{in}\q {\rm int}\, \Om= (0,L)\times (0,1),\\
-\p_x \psi(0,r)=\p_x\psi(L, r)= 0\q &\forall r\in [0,1],\\
\psi(x,0)=\psi(x, 1)=0\q &\forall x\in [0, L].
\end{array}\ee
{Since $-\p_r \xi- \mathbbm{f}_2$ is $C^{\alp}$ in $\Om$ due to \eqref{536}, it follows from the maximum principle and Hopf's lemma that the boundary value problem \eqref{540} has a unique classical solution. Furthermore, the standard Schauder estimate indicates that the classical solution, which is $\psi$, is $C^{2,\alp}$ up to the boundary of $\Om$. Then we obtain from \eqref{536} that $\psi$ satisfies the estimate \eqref{b512}. Note that \eqref{b512} does not mean that $\psi$ as a function in $\mc{N}$ is $C^2$ up to $r=0$. In fact, $\psi$ is not necessarily $C^2$ in $\mc{N}$ up to $r=0$. 
In the next step, we show that $\psi$ satisfies \eqref{limits-psi}, and that the vector field $\psi{\bf e}_{\theta}$ is $C^2$ in $\mc{N}$ up to $r=0$ so that Remark \ref{remark-psi-regularity} implies that ${\bf V}=\psi{\bf e}_{\theta}$ is the unique $C^2$ solution to \eqref{51-vec}.

\smallskip

{\textbf{3.}} In this step, $\psi$ is regarded as a function of the cylindrical coordinates in $\mc{N}$ . By L'Hospital's rule, we have
\be\no
\lim_{r\to 0+}\b(\f{1}r\p_r-\f1{r^2}\b)\psi= \lim_{r\to 0+}\f{r\p_r \psi-\psi}{r^2}=\lim_{r\to 0+}\f{r\p_r^2\psi+\p_r\psi-\p_r\psi}{2r}=\frac{1}{2}\p_r^2\psi(x,0),
\ee
taking the limit $r\to 0+$ to the equation $-\b(\p_{xx}+\frac 1r\p_r(r\p_r) -\f1{r^2}\b) \,\psi= \mathbbm{f}_2(x,r)$, we obtain that
\begin{equation*}
-\lim_{r\to 0+}\frac 32 \p_r^2\psi(x,r)=\mathbbm{f}_2(x,0)+\p_x^2\psi(x,0)= 0\quad\tx{for all}\,\,x\in[0,L],
\end{equation*}
where the second equality is obtained from \eqref{b511} and the condition $\psi(x,0)=0$ on $[0,L]$. This indicates that
\be\no
\p_r^2\psi(x,0)\equiv 0.
\ee
Then a straight forward computation with using $\p_r^2\psi(x,0)\equiv 0$ shows that the vector field ${\bf V}=\psi {\bf e}_{\theta}$ is $C^2$ in $\mc{N}$ with
\begin{equation*}
D^k_{{\bf x}}{\bf V}(x, 0,0)\equiv 0\quad\tx{for}\quad k=0,1,2.
\end{equation*}
By the uniqueness of a $C^2$ solution to \eqref{51-vec}, we finally conclude that \eqref{estimate-V} holds.  
}

\medskip

{\emph{(Method II)}} As another approach to prove the proposition, we modify the arguments used in Lemma 2.2 of \cite{kpr15as}.

 {\textbf{1.}} 
 Let ${\bf V}=V_1({\bf x}){\bf e}_1+V_2({\bf x}){\bf e}_2+V_3({\bf x}){\bf e}_3$ be a $C^2$ solution to \eqref{51-vec}. Here, each ${\bf e}_j$ for $j=1,2,3$ denotes the unit vector in the positive direction of $x_j$-axis for ${\bf x}=(x_1,x_2,x_3)\in \ol{\mc{N}}$. Since ${\bf e}_r\cdot{\bf e}_1=0$, we have $\Delta V_1=0$ in $\mc{N}$, $\der_{x_1} V_1=0$ on $\Ga_0\cup \Ga_L$, and $V_1=0$ on $\Ga_w$. And, this implies that $V_1\equiv 0$ in $\mc{N}$. So it suffices to consider the cylindrical representation of the vector field $V_2{\bf e}_2+V_3{\bf e}_3$ in $\mc{N}$.
\smallskip

Let $\mb{T}$ be a one dimensional torus with period $2\pi$.
As functions of $(x, r,\theta)\in  [0,L]\times [0,1]\times \mb{T}(=:{D}_{\rm cyl})$, we define 
\be\no\begin{cases}
U_r(x,r,\theta):={\bf V}\cdot{\bf e}_r= V_2(x,r\cos\theta, r\sin\theta) \cos\theta + V_3(x,r\cos\theta, r\sin\theta) \sin\theta\\
U_{\theta}(x,r,\theta):={\bf V}\cdot{\bf e}_{\theta}= -V_2(x,r\cos\theta, r\sin\theta) \sin\theta + V_3(x,r\cos\theta, r\sin\theta) \cos\theta.
\end{cases}
\ee 
By the boundary conditions for ${\bf V}$ in \eqref{51-vec}, $(U_r, U_{\theta})$ satisfy
\be\label{BCs-Ucyl}
\der_x(U_r, U_{\theta})=0\quad\tx{on}\,\,\Ga_0\cup\Ga_L,\qquad
(U_r, U_{\theta})=0\quad\tx{on}\,\,\Ga_w.
\ee
Due to $C^{2,\alp}$ regularity of ${\bf V}$ in $\mc{N}$, the functions $U_x$, $U_r$ and $U_{\theta}$ are $C^{2,\alp}$ with respect to the cylindrical variables $(x,r,\theta)$ in ${D}_{\rm cyl}$, and there exists a constant $C$ depending only on $(L,\alp)$ such that
\be\label{estimate-U-cyl-vectors}
\|(U_r, U_{\theta})\|_{2,\alp, D_{\rm cyl}}\le C\|{\bf V}\|_{2,\alp, \mc{N}}
\ee
 Furthermore, $U_x$, $U_r$ and $U_{\theta}$  satisfy
\be\lab{a2}\begin{cases}
\mcl{L}^{\rm cyl}_1(U_r, U_{\theta}):=\b(\p_x^2+\p_r^2+\f1r\p_r+\f1{r^2}\p_{\theta}^2-\f1{r^2}\b) U_r-\f2{r^2}\p_{\theta} U_{\theta}=0,\\
\mcl{L}^{\rm cyl}_2(U_r, U_{\theta}):=\b(\p_x^2+\p_r^2+\f1r\p_r+\f1{r^2}\p_{\theta}^2-\f1{r^2}\b) U_{\theta}+ \f2{r^2}\p_{\theta} U_r=-\mathbbm{f}_2(x,r),
\end{cases}\quad\tx{in}\quad {\rm int}\,D_{\rm cyl}.
\ee
The left-hand sides of the expressions in \eqref{a2} are well-defined  for $r>0$, and are well defined up to $r=0$ by continuation with taking the limits as $r$ tends to $0+$. In taking the limits, note that the fact of $\mathbbm{f}_2(x,0)$ is also essential. 
\smallskip

{\textbf{2.}} For each $n\in \mathbb{N}$,  define functions $U^n_r, U^n_{\theta}$ by
\be\no
U^n(x, r, \theta):=\f{1}{2^n} \sum_{k=0}^{2^n-1} U_r\b(x,r,\theta+\f{2\pi k}{2^n}\b)\,\,,\quad
U^n_{\theta}(x,r,\theta):=\f{1}{2^n} \sum_{k=0}^{2^n-1} U_{\theta}\b(x,r,\theta+\f{2\pi k}{2^n}\b)\quad\tx{in}\,\,D_{\rm cyl}.
\ee
Each $(U^n_r, U^n_{\theta})$ satisfies the following properties:
\begin{itemize}
\item[(i)]  Since the coefficient of each differential operator is independent of $\theta$, it follows from \eqref{a2} that $(U^n_r, U^n_{\theta})$  satisfy
\be\no
(\mcl{L}^{\rm cyl}_1, \mcl{L}^{\rm cyl}_2) (U_r^n, U^n_{\theta})(x,r,\theta)=(0, -\mathbbm{f}_2(x,r))\quad
\tx{in}\,\,D_{\rm cyl};
\ee
\item[(ii)] \eqref{estimate-U-cyl-vectors} yields the estimate
\be\label{estimate-Uns}
\|(U^n_r, U^n_{\theta})\|_{2,\alp, D_{\rm cyl}}\le C\|{\bf V}\|_{2,\alp, \mc{N}};
\ee

\item[(iii)] By definition, each $(U_r^n, U_{\theta}^n)$ satisfies
\be\lab{a7}\begin{cases}
U_r^n(x,r,\theta)= U_r^n(x,r,\theta+\f{2\pi j}{2^n})\\ 
U_{\theta}^n(x,r,\theta)= U_{\theta}^n(x,r,\theta+\f{2\pi j}{2^n}),
\end{cases}\q \forall 0\leq j\leq 2^n-1, n\geq 1.
\ee

\end{itemize}
By \eqref{estimate-Uns} and Arzel\'{a}-ascoli theorem, there exists a sequence $\{n_k\}$ with $\lim_{k\to \infty} n_k=\infty$ such that $\{(U^{n_k}_r, U^{n_k}_{\theta})\}$ converges to functions $(\ti{U}_r, \ti{U}_{\th})$ in $C^{2,\alp/2}(D_{\rm cyl})$. By properties (i) and (ii), $(\ti{U}_r, \ti{U}_{\th})$ satisfy
\begin{equation}
\label{equation-tiU}
(\mcl{L}^{\rm cyl}_1, \mcl{L}^{\rm cyl}_2) (\ti{U}_r, \ti{U}_{\theta})(x,r,\theta)=(0, -\mathbbm{f}_2(x,r))\quad
\tx{in}\,\,D_{\rm cyl},
\end{equation}
and
\begin{equation*}
(\ti{U}_r,\ti{U}_{\theta})(x,r,\theta)= (\ti{U}_r,\ti{U}_{\theta})\b(x,r,\theta+\f{2\pi j}{2^{n_k}}\b)\q\forall 0\leq j\leq 2^{n_k}-1,\,\,\forall k\geq 1.
\end{equation*}
Since $(\ti{U}_r, \ti{U}_{\theta})$  are continuous in $\theta\in \mathbb{T}$, we conclude that $\ti{U}_{r,\theta}(x,r,\theta)=\ti{U}_{r,\theta}(x,r,\theta+2\kappa\pi)$ for any $0\leq \kappa<1$, i.e., $\ti{U}_r, \ti{U}_{\theta}$ are independent of $\theta$. Then the system \eqref{equation-tiU} for $(\ti U_r, \ti U_{\theta})$ are decomposed into two separate elliptic equations:
\be\label{equation-limitU}
\b(\p_x^2+\p_r^2+\f1r\p_r-\f1{r^2}\b) (\ti{U}_r, \ti U_{\theta})(x,r)=(0, -\mathbbm{f}_2(x,r))\quad\tx{in}\,\,D_{\rm cyl}.
\ee 
\smallskip

{\textbf{3.}} To simplify notation, set $\psi(x,r):=\ti{U}_{\theta}(x,r)$ for $(x,r)\in [0,L]\times[0,1]$. Define a vector field ${\bf W}: D_{\rm cyl}\rightarrow \R^3$ by
\begin{equation*}
{\bf W}(x,r,\theta)=\psi(x,r){\bf e}_{\theta}.
\end{equation*}
Note that $\psi$ satisfies the estimate $\|\psi\|_{2,\alp, D_{\rm cyl}}\le C\|{\bf V}\|_{2,\alp,\mc{N}}$ by \eqref{estimate-Uns}.  By \eqref{equation-limitU}, $\psi(x,r)$ can be represented as 
\be\no
\psi(x,r)=r^2\left(\der_x^2+\der_r^2+\frac 1r\der_r\right)\psi(x,r)+r^2\mathbbm{f}_2(x,r)
\ee
for each $r>0$, and the representation is well defined up to $r=0$ by taking limit $r\to 0+$. Furthermore, we obtain that
\be\label{psi-zero-eqn}
\psi(x,0)\equiv 0\quad\tx{for all}\,\,x\in[0,L].
\ee
By repeating Step 3 of ({\emph{Method 1}}), we obtain from \eqref{equation-limitU} and \eqref{psi-zero-eqn} that $\der_r^2\psi(x,0)\equiv 0$ for all $x\in[0,L]$. This implies that the vector field ${\bf W}=\psi{\bf e}_{\theta}$ is in fact $C^2$ in $\mc{N}$ as a vector field of three dimensional Euclidean coordinates. Furthermore, ${\bf W}$ solves \eqref{51-vec}. By the uniqueness of a solution to \eqref{51-vec}, we conclude that 
\begin{equation*}
{\bf V}\equiv {\bf W}\quad\tx{in}\,\,\mc{N}.
\end{equation*}

\epf

\begin{remark}
\label{remark-tvecs}
For $\psi$ in {\emph{Proposition \ref{bl510}}},  the vector field ${\bf t}_1(r,\psi,D\psi):\ol{\mc{N}}\rightarrow \mathbb{R}^3$ given by \eqref{definition-tvecs} is a $C^1$ axisymmetric vector field . Furthermore, there exists a constant $C>0$ depending only on $(L,\al)$ such that
\begin{equation}
\label{estimate-tvec1}
\|{\bf t}_1(r,\psi,D\psi)\|_{1,\al,\mc{N}}\le C\|\psi\|_{2,\al, \Om}.
\end{equation}

\end{remark}
\medskip

\subsection{Transport equation with a div-free vector field }
Finally, we need to solve a linearized version of the problem (\ref{b252}). We regard \eqref{b252} as a problem defined in a two dimensional rectangular domain $\Om=[0,L]\times [0,1]$.

\bp\lab{bl33}
{\it Suppose that a vector field ${\bf M}(x,r)=(M_x(x,r), M_r(x,r))$ satisfies
\be\lab{b349}
\p_x(r M_x)+ \p_r(r M_r)=0\q \forall (x,r)\in \Om=[0,L]\times [0,1],\q M_r(x,0)=M_r(x,1)=0,
\ee
and that $M_x$ satisfies the estimate
\be\lab{b351}
\|M_x\|_{1,\al,\Om}^{(-\al,\{r=1\})}\leq K_0,
\ee
for a constant $K_0>0$. In addition, assume that there exists a constant $\nu^*>0$ satisfying
\be\lab{b352}
M_x\geq \nu^*\q \text{in}\q \Om.
\ee
Then there exists a constant $\e_0>0$ small depending on $(K_0, L)$ such that if $M_r$ satisfies
\be\label{estimate-Mr}
\|M_r\|_{1,\al,\Om}^{(-\al,\{r=1\})}\leq \e_0,
\ee
then the problem
\be\lab{b350}
(M_x\p_x+ M_r\p_r)\mc{W}(x,r)=0\q \text{in}\q \Om,\q \mc{W}=\mc{W}_{en}\q \text{on}\q \Ga_{en}=\partial\Om\cap\{x=0\}
\ee
has a unique solution $\mc{W}\in C^{1,\al}({\Om})$ satisfying
\be\lab{b353}
\|\mc{W}\|_{1,\al,\Om}\leq C^* \|\mc{W}_{en}\|_{1,\al,\Ga_{en}},
\ee
where the constant $C^*$ depends only on $(L,\nu^*, K_0, \e_0, \al)$.
}\ep

\bpf
Set
\be\lab{b354}
w(x, r)\co \int_0^{r} s M_x(x, s) ds\q \text{in}\q \Om.
\ee
It follows from (\ref{b352})--(\ref{b350}) that $w$ satisfies
\be\lab{b355}
\begin{split}
&\begin{cases}
\p_r w= r M_x\geq \nu^* r\\
\p_x w(x,r)= -r M_r(x,r)
\end{cases}
\q \text{in}\q\Om,\\
&\p_x w=0 \q \q\text{on}\q \p \Om\cap\{r=0, 1\},
\end{split}
\ee
and
\be\lab{b356}
\|w\|_{2,\al,\Om}^{(-1-\al,\{r=1\})}\leq C K_0
\ee
where the constant $C$ depending only on $L$.

Since $\p_x w(x,0)= \p_x w(x,1)=0$, we have $w(x,0)=w(0,0)$ and $w(x,1)=w(0,1)$. Also $\p_r w(x, r)= r M_x(x, r)\geq \nu^* r$, then $w(x,r)$ is strictly increasing in $r\in [0,1]$ for each fixed $x\in [0, L]$. This implies that for each $x\in[0,L]$, the closed interval $[w(x,0), w(x, 1)]$ is simply fixed as $[w(0,0), w(0,1)]$.
Therefore, one can uniquely define a function $\vartheta: \Om\rightarrow [0,1]$ to satisfy
\be\lab{b357}
w(x, r)= w(0,\vartheta).
\ee
Set $\mc{G}(r):= w(0,r)$. Since $\mc{G}:[0,1]\rightarrow  [w(0,0), w(0,1)]$ is invertible, the function $\vartheta$ is represented as
\be\lab{b359}
\vartheta(x, r)= \mc{G}^{-1}\circ w(x, r).
\ee
For such a function $\vartheta$, $\mc{W}$ given by
\be\lab{b360}
\mc{W}(x,r)= \mc{W}_{en}(\vartheta(x, r))
\ee
solves \eqref{b350} provided that $\vartheta$ is $C^1$ in $\Om$.
It follows from \eqref{b355} and \eqref{b357} that
\be\lab{b361}
D_{(x,r)} \vartheta(x, r)= \f{D_{(x,r)} w(x, r)}{\p_r w(0,\vartheta(x,r))}= \f{r(-M_r, M_x)(x,r)}{\vartheta(x,r) M_x(0,\vartheta(x, r))}
\ee
unless $\vartheta(x,r)=0$.
By the method of characteristics, $\vartheta(x,r)$ is represented as
\be\lab{b365'}
\vartheta(x,r)-r= -\int_0^x \b(\f{M_r}{M_x}\b)(s,\ka(s;x,r)) ds,
\ee
where $(s, \ka(s; x,r))$ solves
\be\lab{b363'}\begin{cases}
\f{d }{d s}\ka(s; x,r)= \b(\f{M_r}{M_x}\b)(s, \ka(s; x,r))\q \text{for}\,\,0\le s<x\\
\ka(x; x,r)= r.
\end{cases}\ee
for each $(x,r)\in [0,L]\times [0,1]$. Note that $\vartheta(x,r)= \ka(0;x,r)$.
By \eqref{b365'}, one can choose $\e_0>0$ small depending only on $(K_0, L)$ so that if $ r\ge \f 12$, then $\vartheta(x,r)\geq \f 14$ holds, and this implies that
\be\lab{b365}
\|\vartheta\|_{1,\al, [0, L]\times [1/2, 1]}\leq C\|{\bf M}\|_{\al,\Om}.
\ee
To achieve $C^{1,\al}$ estimate of $\vartheta(x,r)$ on $[0, L]\times [0, 1/2]$, we differentiate \eqref{b363'} with respect to $(x,r)$
to get
\be\lab{b363}\begin{cases}
\f{d}{d s}\p_x \ka(s; x,r)= \p_r\b(\f{M_r}{M_x}\b)(s,\ka(s;x,r)) \p_x \ka(s; x,r),\\
\f{d}{d s}\p_r \ka(s; x,r) = \p_r\b(\f{M_r}{M_x}\b)(s,\ka(s;x,r)) \p_r \ka(s; x,r),\\
(\p_x \ka, \p_r\ka)(x; x,r)=\b(-\b(\f{M_r}{M_x}\b)(x,r), 1\b).
\end{cases}\ee
Then we apply Gronwall's  inequality to obtain that
\be\lab{b366}
\|\vartheta\|_{1,\al, [0,L]\times [0,1/2]}\leq C\|{\bf M}\|_{\al,\Om}.
\ee
Finally \eqref{b353} is obtained from combining \eqref{b365} and \eqref{b366} with \eqref{b360}. The uniqueness of a solution can be directly checked by the method of characteristics.

\epf

\begin{remark}
\label{remark-transp}
By \eqref{b349} and \eqref{b360}, we have
\be\no
\partial_r\mc{W}(x,0)=\mc{W}'_{en}(\vartheta(x,0))\partial_r\vartheta(x,0)=\mc{W}'_{en}(0)\partial_r\vartheta(x,0).
\ee
This implies that if $\mc{W}'_{en}(0)=0$, then $\partial_r\mc{W}(x,0)=0$ for all $x\in[0, L]$. Therefore $\mc{W}$ is a $C^1$ axisymmetric function in $\mc{N}$. From this we conclude that if $\mc{W}_{en}$ satisfies the compatibility condition $\mc{W}_{en}'(0)=0$, then the problem  \eqref{b252} with the vector field ${\bf M}=M_x{\bf e}_x+M_r{\bf e}_r$ satisfying \eqref{b253} has a unique $C^1$ axisymmetric solution $\mc{W}\in C^{1,\al}(\ol{\mc{N}})$ with satisfying the estimate
\be\no
\|\mc{W}\|_{1,\al,\mc{N}}\le C^*\|\mc{W}_{en}\|_{1,\al, \partial \mc{N}\cap\{x=0\}},
\ee
for the constant $C^*>0$ depending only on $(L, \nu^*, K_0, \epsilon_0, \al)$.
\end{remark}

\section{Proof of the main theorems}\lab{iteration}
In this section, we first prove Theorem \ref{main1}, then prove Theorem \ref{main}.
\subsection{Proof of Theorem \ref{main1}}

\subsubsection{Step 1:Iteration sets}
Fix $\al\in(0,1)$.

(i) Iteration set for $(S, \mc{K}, \La)$: For a constant $\delta_1>0$ to be determined later, we define
\be\lab{b31}
\mc{P}(\delta_1):=\mc{P}_{\rm pot}(\delta_1)\times \mc{P}_{\rm vort}(\delta_1)
\ee
for
\be\no
\begin{split}
&\mc{P}_{\rm pot}(\delta_1):=\{{\bm \eta}=(S,\mc{K})(x,r)\in [C^{1,\al/2}(\ol{\mc{N}})]^2:  \|(S-S_{0}, \mc{K}-\mc{K}_{0})\|_{1,\al,\mc{N}}\leq \delta_1\},\\
&\mc{P}_{\rm vort}(\delta_1):=\{\La=r\mc{V}(x,r)\in C^{1,\al/2}(\ol{\mc{N}}): \|\mc{V}\|_{1,\al,\Om}\le \delta_1,\,\,
\mc{V}(x,0)=0\,\,\text{for all $x\in[0, L]$} \}.
\end{split}
\ee

(ii) Iteration set for $(\Psi, \phi, \psi)$: For two constants $\delta_2, \delta_3>0$ to be determined later, we define
\be\label{definition-iterset2}
\mc{I}(\delta_2, \delta_3):= & \mc{I}_{\rm pot}(\delta_2)\times \mc{I}_{\rm vort}(\delta_3)
\ee
for
\be\no
\begin{split}
&\mc{I}_{\rm pot}(\delta_2):=\{(\Psi, \phi)(x,r)\in [C_{(-1-\al,\Ga_w)}^{2,\al}(\mc{N})]^2: \|\Psi\|_{2,\al,\mc{N}}^{(-1-\al,\Ga_w)}+\|\phi\|_{2,\al,\mc{N}}^{(-1-\al,\Ga_w)}\leq \delta_2\},\\
&\mc{I}_{\rm vort}(\delta_3):=\{\psi(x,r)\in C^{2,\al}(\ol{\Om}): \|\psi\|_{2,\al, \Om}\le \delta_3, \,\,\psi(x,0)=\psi_{rr}(x,0)=0\,\,\text{for all $x\in[0, L]$}\}.
\end{split}
\ee
By an argument similar to Remark \ref{remark-tvecs}, the following lemma is obtained.
\begin{lemma}
\label{lemma-iter-tvecs}
For each $(\La, \psi)\in \mc{P}_{\rm vort}(\delta_1) \times \mc{I}_{\rm vort}(\delta_3)$, let ${\bf t}_1(r, \psi, D\psi)$ and ${\bf t}_2(r,\La)$ be given by \eqref{definition-tvecs}. Then there exists a constant $C>0$ depending only on $(L, \al)$ such that
\be\no
\|{\bf t}_1(r, \psi, D \psi)\|_{1,\al, \mc{N}}\le C\|\psi\|_{2,\al,\Om},\quad
\|{\bf t}_2(r,\La)\|_{1,\al, \mc{N}}\le C\|\La\|_{1,\al,\mc{N}}.
\ee
\end{lemma}
A direct computation with using Lemma \ref{lemma-iter-tvecs} yields the following result.
\begin{lemma}
\label{lemma-estimate-nonh}
For each $({\bm\eta},\La)\in \mc{P}(\delta_1)$ and $(\Psi ,\phi, \psi)\in \mc{I}(\delta_2, \delta_3)$, let $({\bf F}, f_1)({\rm x}, Q)$, $g(r,\mc{K}-\mc{K}_0, D\phi, {\bf t}(r,\psi, D\psi, \La))$ and $f_2({\rm x}, Q, \p_r{\bm\eta}, \p_r\La)$ be given by \eqref{b25}, \eqref{b26}, \eqref{b28} and \eqref{b211}, respectively, with $Q=({\bm\eta}-{\bm\eta}_0, \Psi, D\phi, {\bf t}(r,\psi, D\psi, \La))$ for ${\bm \eta}_0=(S_0, \mc{K}_0)$. Then there exists a constant $\varepsilon_1>0$ small depending only on the data so that if $\delta_1+\delta_2+\delta_3\le \varepsilon_1$, then we have
\begin{equation}
\label{estimate-nonhs}
\begin{split}
&\|{\bf F}\|_{1,\al, \mc{N}}^{(-\al, \Ga_w)}\le C(\delta_1+\delta_3+\delta_2^2),\\
&\|f_1\|_{\al,\mc{N}}\le C(\delta_1+\delta_3+\om_1(b)+\delta_2^2),\\
&\|g\|_{1,\al,\Ga_L}^{(-\al, \p\Ga_L)}\le C(\delta_1+\delta_3+\om_3(\Phi_{bd},p_{ex})+\delta_2^2),\\
&\|f_2\|_{\al,\mc{N}}\le C\delta_1
\end{split}
\end{equation}
for $\om_1(b)$ and $\om_3(\Phi_{bd}, p_{ex})$ given by \eqref{b213}, where the estimate constant $C$ depends only on the data and $\al$. In addition, $f_2$ satisfies the compatibility condition
\be\no
f_2(x,0)=0\quad\text{for all $x\in[0, L]$}.
\ee
Fix $({\bm\eta}, \La)\in\mc{P}(\delta_1)$. For each $j=1,2$, let ${\bf F}^{(j)}$, $f_1^{(j)}$, $g^{(j)}$ and $f_2^{(j)}$ be defined as above for the fixed $({\bm\eta}, \La)$, and for a fixed $(\Psi^{(j)}, \phi^{(j)}, \psi^{(j)})\in  \mc{I}(\delta_2, \delta_3)$. There exists a constant $\varepsilon_2\in(0,\varepsilon_1]$ depending only on the data and $\al$ so that if $\delta_1+\delta_2+\delta_3\le \varepsilon_2,$ then we have
\begin{equation}
\label{estimate-contraction-add}
\begin{split}
&\|{\bf F}^{(1)}-{\bf F}^{(2)}\|_{\al,\mc{N}}
+\|f_1^{(1)}-f_1^{(2)}\|_{\al,\mc{N}}
+\|g^{(1)}-g^{(2)}\|_{\al,\Ga_L}\\
&\le C\left(\|\psi^{(1)}-\psi^{(2)}\|_{1,\al,\Om}
+(\delta_1+\delta_2+\delta_3)\|(\Psi^{(1)}, \phi^{(1)})-(\Psi^{(2)}, \phi^{(2)})\|_{1,\al, \mc{N}}\right),
\end{split}
\end{equation}
\begin{equation}
\label{estimate-contraction}
\begin{split}
&\|{\bf F}^{(1)}-{\bf F}^{(2)}\|_{1,\al,\mc{N}}^{(-\al, \Ga_w)}
+\|f_1^{(1)}-f_1^{(2)}\|_{\al,\mc{N}}
+\|g^{(1)}-g^{(2)}\|_{1,\al,\Ga_L}^{(-\al, \p\Ga_L)}\\
&\le C\left(\|\psi^{(1)}-\psi^{(2)}\|_{2,\al,\Om}
+(\delta_1+\delta_2+\delta_3)\|(\Psi^{(1)}, \phi^{(1)})-(\Psi^{(2)}, \phi^{(2)})\|_{2,\al, \mc{N}}^{(-1-\al, \Ga_w)}\right),
\end{split}
\end{equation}
and
\begin{equation}
\label{estimate-contraction2}
\|f_2^{(1)}-f_2^{(2)}\|_{\al, \mc{N}}\le C\delta_1(\|\psi^{(1)}-\psi^{(2)}\|_{2,\al,\Om}+\|(\Psi^{(1)}, \phi^{(1)})-(\Psi^{(2)}, \phi^{(2)})\|_{2,\al, \mc{N}}^{(-1-\al, \Ga_w)})
\end{equation}
for a constant $C>0$ depending only on the data and $\al$.
\end{lemma}
\medskip

For a fixed $({\bm\eta}^*, \La^*)\in \mc{P}(\delta_1)$, set
\be\no
Q^*:=({\bm \eta}^*-{\bm\eta}_0, \Psi, \na\phi, {\bf t}(r, \psi, D\psi, \La^*))
\ee
we first solve the following nonlinear boundary value problem for $(\Psi, \phi, \psi)$:
\be\lab{b311}\begin{cases}
L_1(\Psi, \phi)= \text{div}\, {\bf F}({\bf x}, Q^*) \\
L_2(\Psi, \phi)= f_1({\bf x},Q^*)
\end{cases}\q \textit{in}\q \mc{N},\ee
\be\lab{b312}
-\Delta (\psi {\bf e}_{\theta})= f_2({\bf x},Q^*, \p_r(\bm\eta^*, \La^*)) {\bf e}_{\th}\q \textit{in}\q \mc{N}
\ee
with boundary conditions (\ref{b137})-(\ref{b138'}), (\ref{b212}), (\ref{b213}) and
\be\lab{b313}
\p_x \phi= g(r, \mc{K}^*-\mc{K}_0, \na \phi, {\bf t}(r, \psi, D\psi, \La^*))\q \textit{on}\q \Ga_L.
\ee

\subsubsection{Step 2: Well-posedness of the nonlinear boundary value problem for $(\Psi, \phi, \psi)$}

\begin{lemma}\lab{bp32}
{\it
Let $\om_1(b)$, $\om_2(S_{en}, \mc{B}_{en}, v_{en})$, $\om_3(\Phi_{bd}, p_{ex})$ and $\sigma$ be given by \eqref{b213}. Then, there exists a constant $\sigma_5>0$ depending on the data and $\alp$ so that if
\be\lab{b321}
\sigma\leq \si_5,
\ee
then the boundary value problem (\ref{b311})-(\ref{b312}) with boundary conditions (\ref{b137})-(\ref{b138'}), (\ref{b212}), (\ref{b213}) and (\ref{b313}) has a unique solution $(\Psi, \phi, \psi {\bf e}_{\th})\in [C^{2,\al}_{(-1-\al, \Ga_w)}(\mc{N})]^2\times C^{2,\alp}(\ol{\mc{N}}, \R^3)$ with satisfying
\be\lab{b322}
\|\Psi\|_{2,\al,\mc{N}}^{(-1-\al, \Ga_w)}+\|\phi\|_{2,\al,\mc{N}}^{(-1-\al, \Ga_w)}+\|\psi {\bf e}_{\theta}\|_{2,\al,\mc{N}}\le C\sigma
\ee
where the constant $C$ depends only on the data and $\al$. 
}

\begin{proof}
For a fixed $(\tPsi, \tphi, \tpsi)\in \mc{I}(\delta_2, \delta_3)$, we set
\be\no
\tilde{Q}^*:=({\bm\eta}^*-{\bm\eta}_0, \tPsi, \na\tphi, {\bf t}(r,\tpsi, D\tpsi, \La^*)),
\ee
and solve the following associated linear boundary value problem
\be\label{lbvp-1}
\begin{split}
&\begin{cases}
L_1(\Psi, \phi)= \text{div}\, {\bf F}({\bf x}, \tilde{Q}^*) \\
L_2(\Psi, \phi)= f_1({\bf x},\tilde{Q}^*)\\
-\Delta(\psi{\bf e}_{\theta})= f_2({\bf x},\tilde{Q}^*, \p_r(\bm\eta^*, \La^*)){\bf e}_{\theta}
\end{cases}\q \textit{in}\q \mc{N},\\
&\p_x \phi= g(r, \mc{K}^*-\mc{K}_0, \na \tphi, {\bf t}(r, \psi, D \tpsi, \La^*))\q \textit{on}\q \Ga_L,\\
&\text{with boundary conditions (\ref{b137})-(\ref{b138'}), (\ref{b212}), (\ref{b213-new}). }
\end{split}
\ee

By \eqref{b213}, Lemma \ref{b142}, Proposition \ref{bl510} and Lemma \ref{lemma-estimate-nonh}, the linear boundary value problem \eqref{lbvp-1} has a unique solution $(\Psi, \phi, \psi{\bf e}_{\theta})\in [C^{2,\al}_{(-1-\al, \Ga_w)}(\mc{N})]^2\times C^{2,\alp}(\ol{\mc{N}}, \R^3)$ with satisfying the estimates
\begin{equation}
\label{estimate-sol-lbvp1}
\begin{split}
&\|\Psi\|_{2,\al,\mc{N}}^{(-1-\al, \Ga_w)}+\|\phi\|_{2,\al,\mc{N}}^{(-1-\al, \Ga_w)}\le \mc{C}_1(\delta_1+\delta_3+\delta_2^2+\sigma),\\
&\|\psi {\bf e}_{\theta}\|_{2,\al,\mc{N}}\le \mc{C}_1\delta_1
\end{split}
\end{equation}
for a constant $\mc{C}_1>0$ depending only on the data and $\al$.

We choose $\delta_3$ as
\be\label{choice-d3}
\delta_3=2\mc{C}_1\delta_1.
\ee
Under such a choice of $\delta_3$, if it holds that
\be\label{choice-d-cond1}
\mc{C}_1\left((1+2\mc{C}_1)\frac{\delta_1}{\delta_2}+\delta_2+\frac{\si}{\delta_2}\right)\le \frac 12,
\ee
then we have $(\Psi, \phi, \psi)\in \mc{I}(\delta_2, \delta_3)$. We define a mapping $\mathfrak{I}^{({\bm\eta}^*, \La^*)}$ by
\be\no
\mathfrak{I}_1^{({\bm\eta}^*, \La^*)}(\tPsi, \tphi, \tpsi):=(\Psi, \phi, \psi)\quad\text{for each $(\tPsi, \tphi, \tpsi)\in\mc{I}(\delta_2, \delta_3)$},
\ee
where $(\Psi, \phi, \psi{\bf e}_{\th})\in C^{2,\alp}_{(-1-\alp, \Ga_w)}(\mc{N})]^2\times C^{2,\alp}(\ol{\mc{N}},\R^3)$ is the unique axisymmetric solution to \eqref{lbvp-1}. Then $\mathfrak{I}_1^{({\bm\eta}^*, \La^*)}$ maps $\mc{I}(\delta_2, \delta_3)$ into itself. We regard $\mc{I}(\delta_2,\delta_3)$ as a compact and convex subset of $[C^{1,\alp/2}(\ol{\mcl{N}})]^2\times C^{2,\alp/2}(\ol{\Om})$. Since $\mathfrak{I}_1^{({\bm \eta}^*, \La^*)}:\mc{I}(\delta_2, \delta_3)\rightarrow \mc{I}(\delta_2, \delta_3)$ is continuous in $[C^{1,\alp/2}(\ol{\mc{N}})]^2\times C^{2,\alp/2}(\ol{\Om})$, the Schauder fixed point theorem implies that $\mathfrak{I}_1^{({\bm \eta}^*, \La^*)}$ has a fixed point $(\Psi, \phi, \psi)\in\mc{I}(\delta_2, \delta_3)$.

Let $(\Psi^{(1)}, \phi^{(1)}, \psi^{(1)})$ and $(\Psi^{(2)}, \phi^{(2)}, \psi^{(2)})$ be two fixed points of $\mathfrak{I}_1^{({\bm \eta}^*, \La^*)}$. Then it follows from Lemma \ref{bl42}, Proposition \ref{bl510}, Lemma \ref{lemma-estimate-nonh}  and \eqref{choice-d3} that
\begin{equation*}
\|\psi^{(1)}-\psi^{(2)}\|_{2,\alp,\Om}
\le \mcl{C}_2\delta_1 (\|\psi^{(1)}-\psi^{(2)}\|_{2,\al,\Om}+\|(\Psi^{(1)}, \phi^{(1)})-(\Psi^{(2)}, \phi^{(2)})\|_{2,\al, \mc{N}}^{(-1-\al, \Ga_w)}),
\end{equation*}
and
\begin{equation*}
\begin{split}
&\|(\Psi^{(1)}, \phi^{(1)})-(\Psi^{(2)}, \phi^{(2)})\|_{2,\al, \mc{N}}^{(-1-\al, \Ga_w)}\\
&\le \mcl{C}_2\left(\|\psi^{(1)}-\psi^{(2)}\|_{2,\al,\Om}
+(\delta_1+\delta_2)\|(\Psi^{(1)}, \phi^{(1)})-(\Psi^{(2)}, \phi^{(2)})\|_{2,\al, \mc{N}}^{(-1-\al, \Ga_w)}\right).
\end{split}
\end{equation*}
for a constant $\mcl{C}_2>0$ depending only on the data and $\alp$. If it holds that
\begin{equation}
\label{choice-d-cond2}
\mcl{C}_2\delta_1\le \frac 12,
\end{equation}
then we obtain from the previous two estimates that
\begin{equation*}
\|(\Psi^{(1)}, \phi^{(1)})-(\Psi^{(2)}, \phi^{(2)})\|_{2,\al, \mc{N}}^{(-1-\al, \Ga_w)}\le \mcl{C}_3(\delta_1+\delta_2) \|(\Psi^{(1)}, \phi^{(1)})-(\Psi^{(2)}, \phi^{(2)})\|_{2,\al, \mc{N}}^{(-1-\al, \Ga_w)}
\end{equation*}
for a constant $\mcl{C}_3>0$ depending only on the data and $\alp$. We conclude that the fixed point is unique provided that
\begin{equation}
\label{choice-d-cond3}
\mcl{C}_3(\delta_1+\delta_2)\le \frac 12.
\end{equation}

We now make choices of $\delta_1$ and $\delta_2$. For $\sigma_5>0$ to be specified later, we choose $(\delta_1,\delta_2 )$ as
\be\label{choice-delta-final}
\delta_1=\frac{\delta_2}{10\mcl{C}_1(1+2\mcl{C}_1)},\quad \delta_2=\frac{10}{\mcl{C}_1}\sigma_5
\ee
for the constant $\mcl{C}_1$ from \eqref{estimate-sol-lbvp1}. Then, \eqref{choice-d-cond1} holds whenever $\sigma\in(0, \sigma_5]$. Finally, we choose $\sigma_5$ as
\be\label{choice-sigma5-final}
\sigma_5 =\mcl{C}_1^2(1+2\mcl{C}_1)\min\{ \frac{1}{2\mcl{C}_2}, \frac{1}{2\mcl{C}_3\left(1+10\mcl{C}_1(1+2\mcl{C}_1)\right)}\},
\ee
so that \eqref{choice-d-cond2} and \eqref{choice-d-cond3} hold. The proof is completed.
\end{proof}
\end{lemma}

\subsubsection{Step 3: Existence of a solution to the problem (\ref{b141})-(\ref{b146}) with (\ref{b138'})-(\ref{b320'}) }

In the proof of Proposition \ref{bp32}, we have shown that there exist constants $\sigma_5>0$ and $\mcl{C}>0$ depending only on the data and $\al$ so that if $\sigma\le \sigma_5$, then for each $\mc{W}^*=({\bm\eta}^*, \La^*)\in\mc{P}(\delta_1)$, there exists a unique axisymmetric weak solution $(\Psi, \phi,\psi)\in[C^{2,\alp}_{(-1-\alp,\Ga_w)}(\mc{N})]^2\times [C^{0,1}(\ol{\mc{N}})\cap C^{2,\alp}(\mc{N})\cap C^{2,\alp}(\ol{\Om})]$ to \eqref{b311}-\eqref{b312} with boundary conditions \eqref{b137}--\eqref{b138'}, \eqref{b212}, \eqref{b213} and \eqref{b313}. For such functions $(\Psi, \phi, \psi)$, define a vector field
\be\lab{b348}
{\bf M}^{\mc{W}^*}:=H(S^*, \mc{K}^*+\Phi_0+ \Psi-\f 12 |{\bf q}(r,\psi, D\psi, D\varphi_0+D\phi, \La^*)|^2)|{\bf q}(r,\psi, D\psi, D\varphi_0+D\phi, \La^*),
\ee
for ${\bf q}$ and $H$ defined by \eqref{def-q} and \eqref{definition-H}, respectively. Here, $D$ denotes $D=(\p_x, \p_r)$. By \eqref{b322} and \eqref{choice-delta-final}, ${\bf M}={\bf M}^{\mc{W}^*}$ satisfies \eqref{b349}, and the estimate
\be\lab{b51}
{{
\|{\bf M}^{\mc{W}^*}-J_0{\bf e}_x\|_{1,\alp,\Om}^{(-\alp, \{r=1\})} \le C(\sigma+\sigma_5)}}
\ee
for a constant $C>0$ depending only on the data and $\alp$, where $J_0$ is the momentum density in $x$-direction of the background solution. One can reduce $\sigma_5$ from \eqref{choice-sigma5-final} depending only on the data and $\alp$ so that whenever $\sigma\le \sigma_5$, it follows from Proposition \ref{bl33} and Remark \ref{remark-transp} that the problem \eqref{b252} with ${\bf M}={\bf M}^{\mc{W}^*}$ has a unique solution $\mc{W}=(S, \mc{K}, \La)\in [C^{1,\alp}(\ol{\mc{N}})]^3$ with satisfying
\be\label{estimate-transps}
\|\mc{W}-(S_0, \mc{K}_0, 0)\|_{1,\alp,\mc{N}}
 \le C_T^{(1)}(\om_2(S_{en}, \mc{B}_{en}, v_{en})+\om_3(\Phi_{bd}), p_{ex})\le C_T^{(1)}\sigma
\ee
for a constant $C_T^{(1)}>0$ depending only on the data and $\alp$.

By \eqref{b320'} and \eqref{b360}, $\La$ is represented as
\be\no
\La(x,r)= \La_{en}(\vartheta(x,r))= \vartheta(x,r) v_{en}(\vartheta(x,r)),
\ee
where $\vartheta(x,r)$ is given by \eqref{b359} associated with ${\bf M}={\bf M}^{\mcl{W}^*}$.
Set $\mc{V}$ as
\be\no
\mc{V}(x,r)= \begin{cases}
\f{\vartheta(x,r)}{r} v_{en}(\vartheta(x,r)), \q&\textit{for}\q (x,r)\in [0,L]\times (0,1],\\
0,\q&\textit{for}\q (x,r)\in \{(x,0): x\in [0,L]\}.
\end{cases}
\ee
By Lemma \ref{bl33}, we have $\vartheta(x,0)\equiv 0$ for all $x\in [0,L]$ and $\|\vartheta\|_{1,\al,\ol{\mc{N}}}\leq C\|{\bf M}^{\mc{W}^*}\|_{\al, {\mc{N}}}$. By using \eqref{compatibility-nu} and the representation
$
\p_r \mc{V}(x,r)= \f{\vartheta}{r} v_{en}'(\vartheta)\p_r \vartheta + \p_r \vartheta\f{v_{en}(\vartheta)}{r}- \f{\vartheta}{r}\cdot \f{v_{en}(\vartheta)}{r},
$ we get $\lim_{r\to 0+} \p_r\mc{V}(x,r)=\vartheta_r^2(x,0)\nu'_{en}(0)$. From this, it can be directly checked that
$\mc{V}\in C^{1,\al}(\ol{\Om})$, and
\be\label{estimate-mcV}
\|\mc{V}\|_{1,\al,\ol{\mc{N}}}\leq C^{(2)}_{T}\si
\ee
for a constant $C_T^{(2)}>0$ depending only on the data and $\alp$. For the constant $\mcl{C}_1$ from \eqref{estimate-sol-lbvp1}, we choose $\sigma_6\in(0,\sigma_5]$ as
\be\label{choice-d6}
\sigma_6=\sigma_5\min\left\{1, \frac{1}{\mcl{C}_1^2(1+2\mcl{C}_1)\max\{C_T^{(1)}, C_T^{(2)}\}}\right\}.
\ee
It follows from \eqref{choice-delta-final}, \eqref{estimate-transps} and \eqref{estimate-mcV} that if $\sigma\le \sigma_6$, then we have $\mc{W}\in \mc{P}(\delta_1)$.

For each $\sigma\in(0, \sigma_6]$, we define an iteration mapping $\mc{J}: \mc{P}(\delta_1)\to  \mc{P}(\delta_1)$ by
\be\lab{b450}
\mc{J} \mc{W}^*= \mc{W},
\ee
where $\mc{W}$ is the solution to the problem \eqref{b252} with ${\bf M}={\bf M}^{\mc{W}^*}$. From using Lemma \ref{bp32} and the facts that $\mc{P}(\delta_1)\times \mc{I}(\delta_2, \delta_3)$ is compact in $[C^{1,\alp/2}(\ol{\mc{N}})]^3\times [C^{2,\alp/2}_{(-1-\alp/2, \Ga_w)}(\mc{N})]^2\times C^{2,\alp/2}(\ol{\Om})$ and that the problem \eqref{b252} with ${\bf M}={\bf M}^{\mc{W}^*}$ has a unique solution for each $\mc{W}^*\in\mc{P}(\delta_1)$, we obtain that $\mc{J}: \mc{P}(\delta_1)\to  \mc{P}(\delta_1)$ is continuous in $[C^{1,\alp/2}(\ol{\mc{N}})]^3$. Since $\mc{P}(\delta_1)$ is convex and compact in
$[C^{1,\alp/2}(\ol{\mc{N}})]^3$, we conclude from Schauder fixed point theorem that $\mc{J}$ has a fixed point $\mc{W}=(S, {\mc{K}}, \La)\in \mc{P}(\delta_1)$. For such $\mc{W}$, let $(\Psi, \phi,\psi)$ be the unique fixed point of $\mathfrak{I}_1^{\mc{W}}$ in $\mc{I}(\delta_2, \delta_3)$. Let us set $(\Phi, \varphi)=(\Phi_0, \varphi_0)+(\Psi, \phi)$. Then, $(S, {\mc{K}}, \La, \Phi, \varphi, \psi)$ solves  the problem (\ref{b141})-(\ref{b146}) with (\ref{b138'})-(\ref{b320'}) provided that $\sigma\le \sigma_6$.

\subsubsection{Step 4: Uniqueness of a solution to the problem (\ref{b141})-(\ref{b146}) with (\ref{b138'})-(\ref{b320'})}
For each $j=1,2$, set
\be\no
\mc{U}^{(j)}:=(\Phi^{(j)}, \phi^{(j)}, \psi^{(j)}),\quad \mc{W}^{(j)}:=(S^{(j)}, \mc{K}^{(j)}, \La^{(j)}).
\ee
Let $(\mc{U}^{(j)}, \mc{W}^{(j)}) (j=1,2)$ be two solutions to the problem (\ref{b141})-(\ref{b146}) with (\ref{b138'})-(\ref{b320'}), and set
\be\no
\begin{split}
&d_{\rm el}:=\|(\Psi^{(1)}-\Psi^{(2)}, \phi^{(1)}-\phi^{(2)})\|_{1,\alp,\mc{N}}+\|\psi^{(1)}-\psi^{(2)}\|_{1,\alp,\Om},\\
&d_{\rm trans}:=\|\mc{W}^{(1)}-\mc{W}^{(2)}\|_{\alp,\mc{N}}
\end{split}
\ee
Assume that $\om_4=\om_4(S_{en}, \mc{B}_{en}, \nu_{en}, \Phi_{bd})$ given by \eqref{e3} is finite.

For each $j=1,2$, let ${\bf F}^{(j)}, f_1^{(j)}, f_2^{(j)}$ and $g^{(j)}$ be given by (\ref{b25}), (\ref{b26}), (\ref{b28}) and (\ref{b211}) with
\be\no
Q= ((S^{(j)}-S_0, \mc{K}^{(j)}-\mc{B}_0),\Psi^{(j)}, \na\phi^{(j)}, {\bf t}(r,\psi^{(j)}, D_{(x,r)}\psi^{(j)}, \La^{(j)})), \p_r \mc{W}^{(j)}).
\ee

It follows from \eqref{b25}, \eqref{b26}, \eqref{b28},
(\ref{b411}) and \eqref{estimate-contraction-add} that
\be\lab{94}
\begin{split}
&\|(\phi^{(1)}-\phi^{(2)}, \Psi^{(1)}-\Psi^{(2)})\|_{1,\alp, \mc{N}}\\
&\leq C\left(
\|\psi^{(1)}-\psi^{(2)}\|_{1,\alp,\Om}+\sigma \bigl(\|(\phi^{(1)}-\phi^{(2)}, \Psi^{(1)}-\Psi^{(2)})\|_{1,\alp, \mc{N}}+d_{\rm trans}\bigr)
\right).
\end{split}
\ee
where the constant $C$ depends only on the data and $\alp$. In the following estimates, each estimate constant $C$ may vary from line to line, but it is regarded to be depending only on the data and $\alp$ unless otherwise specified.

By \eqref{def-xi-psi} and \eqref{536}, we have
\be\lab{95}
\|\psi^{(1)}-\psi^{(2)}\|_{1,\alp,\Om}\leq C \|f_2^{(1)}-f_2^{(2)}\|_{\alp, \mc{N}},
\ee
and a straightforward calculation with using \eqref{b322-1} and \eqref{b211}  yields
\be\lab{96}
&&\begin{array}{ll}
\|f_2^{(1)}-f_2^{(2)}\|_{\alp, \mc{N})}
\leq C \left(
\sigma (d_{\rm el}+d_{\rm trans})
+\|\der_r \mc{W}^{(1)}-\der_r \mc{W}^{(2)}\|_{\alp, \Om}
\right ).
\end{array}\ee
By Lemma \ref{bl33}, we have
\be\label{expression-W}
{ \mc{W}^{(j)}= \left(S_{en}(\vartheta^{(j)}),\, \mc{B}_{en}(\vartheta^{(j)})-\Phi_{bd}(0, \vartheta^{(j)}),\, \vartheta^{(j)}\nu_{en}(\vartheta^{(j)})\right)},
\ee
 where $\vartheta^{(j)}$ is given by \eqref{b359} associated with the vector field
\be\no
{\bf M}^{(j)}:= \left(\na\var_0+\nabla \phi^{(j)}+{\bf t}(r,\psi^{(j)}, D\psi^{(j)}, \La^{(j)})\right) H^{(j)}
\ee
with $H^{(j)}:=H(S^{(j)}, \mc{K}^{(j)}+\Phi_0+\Psi^{(j)}-\f 12 |\na\var_0+\nabla \phi^{(j)}+{\bf t}(r,\psi^{(j)}, D\psi^{(j)}, \La^{(j)})|^2)$ for $H$ given by \eqref{definition-H}. Then we get
\be\lab{99}
\|\mc{W}^{(1)}-\mc{W}^{(2)}\|_{\alp,\mc{N}}
\le C\sigma \|\vartheta^{(1)}-\vartheta^{(2)}\|_{\alp,\Om}.
\ee
Furthermore it can be directly checked from \eqref{b355} and \eqref{b357} that
\be\lab{911}
\|\vartheta^{(1)}-\vartheta^{(2)}\|_{\alp, \mc{N}}\leq C(d_{\rm el}+d_{\rm trans}).
\ee
We differentiate \eqref{expression-W} with respect to $r$, then apply \eqref{b365}, \eqref{b366} and \eqref{911} to get
\be\no
\|\p_r \mc{W}^{(1)}-\p_r\mc{W}^{(2)}\|_{\alp, \Om}\le
C\left(\om_4(S_{en},\mc{B}_{en},\nu_{en}, \Phi_{bd})(d_{\rm el}+d_{\rm trans})
+\sigma \|\p_r\vartheta^{(1)}-\p_r\vartheta^{(2)}\|_{\alp, \Om}
\right).
\ee
Finally, straightforward computations with using \eqref{b322-1}, \eqref{b361}, \eqref{b363} and Gronwall's inequality yields that
\be\lab{918}
\|\p_r (\vartheta^{(1)}-\vartheta^{(2)})\|_{\alp, \Om}\leq C
(d_{\rm el}+d_{\rm trans}).
\ee

From the estimates \eqref{94}--\eqref{918}, it is obtained that
\be\label{estimate-final-cont}
(d_{\rm el}+d_{\rm trans})\le
\mc{C}_*\left(\sigma+\om_4(S_{en},\mc{B}_{en},\nu_{en}, \Phi_{bd})\right)(d_{\rm el}+d_{\rm trans})
\ee
for a constant $\mc{C}_*>0$ depending only on the data and $\alp$. For $\sigma_6$ from \eqref{choice-d6}, if
\be\no
\sigma+\om_4(S_{en},\mc{B}_{en},\nu_{en}, \Phi_{bd})\le \min\{\sigma_6, \frac{4}{5\mcl{C}_*}\},
\ee
then \eqref{estimate-final-cont} implies that $(\mcl{U}^{(1)}, \mc{W}^{(1)})\equiv(\mcl{U}^{(2)}, \mc{W}^{(2)})$ in $\mc{N}$.
The proof of Theorem \ref{main1} is completed by choosing $\sigma_3$ and $\sigma_4$ as
\be\no
{{
\sigma_3=\sigma_6,\qquad \sigma_4=\min\{\sigma_6, \frac{4}{5\mcl{C}_*}\}.}}
\ee

\hfill $\Box$

\subsection{Proof of Theorem \ref{main}}
Let $\sigma_3$ be from Theorem \ref{main1}. We choose $\sigma_1$ from \eqref{b212} as
\be\no
\sigma_1=\sigma_3.
\ee
Given $(b, S_{en}, \mc{B}_{en}, \nu_{en}, \Phi_{bd}, p_{ex})(r)$ satisfying \eqref{b212} and \eqref{compatibility-nu}, let $\Sigma=(\varphi, \psi, \Phi, S, \mc{K}, \La)$ be a solution to the problem
\eqref{b141}--\eqref{b146}  with boundary conditions \eqref{b138'}--\eqref{b320'}.
By Theorem \ref{main1}, such a solution $\Sigma$ exists, and it satisfies the estimate \eqref{b322-1}. We define ${\bf u}$ and $\rho$ by \eqref{b110} and $\eqref{b139}$, respectively. Then $(\rho, {\bf u}, S,\Phi)$ solve Problem \ref{problem}. We particularly emphasize that Remarks \ref{remark-tvecs} and \ref{remark-transp} imply that the vector field ${\bf u}$ given by \eqref{b110} is a $C^1$ axisymmetric vector field in $\mc{N}$. Therefore, \eqref{b322-1} implies that $(\rho, {\bf u}, S, \Phi)$ satisfy the estimate \eqref{e1}. Furthermore, the choice of $\sigma_1=\sigma_3(=\sigma_6)$ ensures that  $u_x:={\bf u}\cdot {\bf e}_x>0$ and $\rho>0$ hold in $\ol{\mc{N}}$.

For each $j=1,2$, let $(\rho^{(j)}, {\bf u}^{(j)}, p^{(j)}, \Phi^{(j)})$ be a solution to Problem \ref{problem} with satisfying the estimate \eqref{e1}.
For each $j=1,2$, we write ${\bf u}^{(j)}$ as
$
{\bf u}^{(j)}=u^{(j)}_x{\bf e}_x+ u^{(j)}_r\cdot{\bf e}_r+ u_{\theta}^{(j)}{\bf e}_{\theta},
$
and solve the following linear boundary value problem:
\be\lab{e11}
-\Delta (\psi^{(j)} {\bf e}_{\theta})=  (\p_x u_r^{(j)}- \p_r u_x^{(j)}) {\bf e}_{\theta}
\ee
with the boundary conditions
\be\lab{e12}
\p_x \psi^{(j)}=0\q \text{on}\q \Ga_0\cup\Ga_L,\q\q\psi^{(j)}=0\q \text{on}\q \Ga_w\cup\{(x,0):0\le x\le L\}.
\ee
By Proposition \ref{bl510}, the boundary value problem \eqref{e11}--\eqref{e12} has a unique axisymmetric solution $\psi^{(j)}{\bf e}_{\th} \in C^{2,\alp}(\ol{\mc{N}}, \R^3)$.
We define functions $\var^{(j)}$ and $\La^{(j)}$ by
\be\lab{e13}
\var^{(j)}(x,r)= \int_0^x u_x^{(j)}(y,r)- \f{1}{r}\p_r(r\psi^{(j)}(y,r))\, dy,\q \La^{(j)}(x,r)= ru_{\theta}(x,r)\q\text{in}\,\,\mc{N}.
\ee
By using (\ref{e11}), one can directly check that
\be\no
 \p_r \var^{(j)}(x,r)=u_r(x,r)+\p_x \psi^{(j)}(x,r)\quad\text{in}\,\,\mc{N}.
\ee
For each $j=1,2$, we also define
\be\no
S^{(j)}= c_v \log \b(\f{p^{(j)}}{A\rho^{(j)}}\b),\q\q \mc{K}^{(j)}= \frac 12|{\bf u}^{(j)}|^2+\frac{\gamma p^{(j)}}{(\gamma-1)\rho^{(j)}}-\Phi^{(j)}\q\text{in}\,\,\mc{N}.
\ee
Then each $(\var^{(j)},\psi^{(j)},\Phi^{(j)}, S^{(j)}, \Phi^{(j)}, \La^{(j)})$ solves the problem (\ref{b141})-(\ref{b146}) with (\ref{b138'})-(\ref{b320'}). Furthermore, it follows from \eqref{e1} that each $(\var^{(j)},\psi^{(j)},\Phi^{(j)}, S^{(j)}, \Phi^{(j)}, \La^{(j)})$ satisfies the estimate \eqref{b322-1}. Finally, we choose $\sigma_2$ from \eqref{e2} as $\sigma_2=\sigma_4$ for $\sigma_4$ from \eqref{b323} so that Theorem \ref{main} implies that
 $\Xi^{(1)}=\Xi^{(2)}$. Hence $(\rho^{(1)},{\bf u}^{(1)}, p^{(1)}, \Phi^{(1)})=(\rho^{(1)},{\bf u}^{(2)}, p^{(2)}, \Phi^{(2)})$ in $\mc{N}$. This completes the proof of Theorem \ref{main}. \hfill $\Box$.

\bigskip

{\bf Acknowledgement.}
The research of Myoungjean Bae was supported in part by Samsung Science and Technology Foundation
under Project Number SSTF-BA1502-02.
The research of Shangkun Weng was supported in part by Priority Research Centers Program through the National Research Foundation of Korea(NRF) funded by the Ministry of Education(2015049582).


\end{document}